\documentclass[a4paper, 11pt]{amsart}
\pdfoutput=1
\usepackage{amssymb,amsmath,txfonts,mathrsfs,titletoc,appendix}
\usepackage{graphicx}
\usepackage{latexsym}
\usepackage[alphabetic]{amsrefs}
\usepackage{geometry}
\usepackage{fancyhdr}

 \newtheorem{thm}{Theorem}[section]
 \newtheorem{cor}[thm]{Corollary}
 \newtheorem{lem}[thm]{Lemma}
 \newtheorem{prop}[thm]{Proposition}
 \newtheorem{defn}[thm]{Definition}
 \newtheorem{rem}[thm]{Remark}
 \numberwithin{equation}{section}
\newtheorem{lem*}{Lemma}
\newtheorem{cor*}{Corollary}
\newenvironment{prooff}{\medskip \noindent
{\bf Proof.}}{\hfill \rule{.5em}{1em}
\\}
\newenvironment{pf}{\medskip \noindent
{\bf Proof of Theorem \ref{mainthm1}.}}{\hfill \rule{.5em}{1em}
\\}
\newenvironment{pf2}{\medskip \noindent
{\bf Proof of Theorem \ref{dimdiffthm}.}}{\hfill \rule{.5em}{1em}
\\}

\newenvironment{pf4}{\medskip \noindent
{\bf Proof of Claim}.}{\hfill \rule{.5em}{1em}
\\}

\begin{document}

\title[To control distributions of mass-minimizing currents]{How to construct metrics to control distributions of homologically mass-minimizing currents}

\author{Yongsheng Zhang}
\address{\em{Current address:} School of Mathematics and Statistics, Northeast Normal University, 5268 RenMin Street, NanGuan District, ChangChun, JiLin 130024, P.R. China }
\email{yongsheng.chang@gmail.com}
\date{\today}
\date{\today}

\begin{abstract}
By  Federer and Fleming \cite{FF} there exist at least one mass-minimizing normal current in every real-valued homology class of a Riemannian manifold.
However the regularity of the mass-minimizing currents and their distributions may generally be quite complicated. 
In this paper we shall study how to construct nice metrics so that (as functionals over smooth forms) 
almost all homologically mass-minimizing currents are just linear combinations over smooth submanifolds.

\end{abstract}


\maketitle
\titlecontents{section}[0em]{}{\hspace{.5em}}{}{\titlerule*[1pc]{.}\contentspage}
\titlecontents{subsection}[1.5em]{}{\hspace{.5em}}{}{\titlerule*[1pc]{.}\contentspage}
{\setcounter{tocdepth}{1} \small \tableofcontents}
\section{Introduction}\label{Section1} 

 A calibration on a Riemannian manifold (see \S\ref{cg}) is a closed differential $m$-form $\phi$
 whose value at every point on every unit $m$-plane is at most one.
 The fundamental theorem of calibrated geometry in \cite{HL2} asserts that 
 a smooth $m$-dimensional oriented submanifold (or more generally an integral or normal current) without boundary
 for which $\phi$ has value one at almost every point on every (generalized) unit tangent plane is mass-minimizing in its homology class
 of the complex of normal currents. Therefore the theory of calibrations is a powerful tool for proving the property of homological mass-minimality.

Notice that besides a background smooth manifold there are two extra slots, namely a balanced pair of a closed form and a metric, in a calibrated manifold.
In this paper we shall create such balanced pairs for various situations.
Particular interests are focused on two kinds of metric constructions based upon any given metric.
One is the horizontal change in \S\ref{H} and the other is the conformal change in \S\ref{ccm}.
Given a homologically nontrivial, oriented, compact, connected smooth submanifold $M$,
we show that one can alter any metric by the horizontal change or the conformal change such that
$M$ becomes homologically mass-minimizing with respect to the resulted metric.
It should be pointed out that Tasaki \cite{Tasaki} first proved the existence of metrics making $M$ homologically mass-minimizing.
Although the difference between Tasaki's result and ours reflects in the metric slot,
our improvement essentially comes from (the construction of potential) calibrations.
Lemma \ref{fperp} provides us with a well-behaved potential local calibration,
relying on which
explicit metric gluing constructions
can be made for global calibrations.
The approach is feasible since submanifolds are ``isolated" vs. a foliation
(cf. \cite{Gluck}, \cite{S3} and \cite{HL1}).

 The case of a constellation of mutually disjoint submanifolds possibly of different dimensions
 is also studied in \S\ref{H} and \S\ref{ccm}.
 However only one calibration is constructed for each dimension.
 In order to control distributions of (almost all) homologically mass-minimizing currents
 we construct a metric under which enough calibrations are guaranteed for each dimension.
 Theorem \ref{1} asserts that if $X^n$ is oriented with betti numbers $b_k<\infty$ for $1\leq k< n/2$,
then there exists some metric in any conformal class such that
all homologically area-minimizing currents of dimension $<n/2$ %$<\frac{n}{2}$
are currents of $\mathbb R$-weighted integration over mutually disjoint finitely many oriented, compact, connected submanifolds.
By combination of the horizontal change and the conformal change,
Theorem \ref{2} confirms the general existence of such metrics allowing $k$ to vary up to $n-3$ (without the restriction to a fixed conformal class).
\\{\ }

%%%%%%%%%%%%%%%
%%%%%%%%%%%%%%%%%%%%%%%%%%%%%%%%%%%%%%%%%%%%%%%%%%%%%%%%%%%%%%%%%%%%%%%%%%
$\text{\sc Acknowledgement}.$
This paper is an expansion of part of the author's Ph.D. thesis at Stony Brook University.
He is deeply indebted to his advisor Professor H. Blaine Lawson, Jr. for suggesting this series of questions and constant encouragement.
He also would like to thank Professor Frank Morgan and Professor Dennis Sullivan for helpful comments.
Part of the work was polished during the author's visit to the MSRI in Fall 2013.
\\
{\ }

\section{Preliminaries}\label{P} % Main section title

\subsection{Calibrated Geometry}\label{cg}
Let us review some fundamental concepts and results that we will need in this paper.
Readers are referred to \cite{HL2} for further understanding of calibrated geometry
and to \cite{FM} for an overview of geometric measure theory.

\begin{defn}\label{comassdef}
Let $\phi$ be a smooth $m$-form on a Riemannian manifold $(X,g)$. At a point $x\in X$ we define the $\mathrm{\mathbf{comass}}$ of $\phi_x$ to be
\begin{equation*}
\|\phi\|_{x,g}^*=\max \ \{\phi_x( \overrightarrow V_x) : \overrightarrow V_x \ \text{is a unit simple m-vector at x}\}
\end{equation*}
where {``simple"} means
$\overrightarrow V_x=e_1\wedge e_2\cdots \wedge e_m$
for some 
$e_i\in T_xX$.
%Furthermore, the $\mathrm{\mathbf{comass}}$ of $\phi$ (over $X$) is defined as
%$$\|\phi\|_{X,g}^*=\sup\|\phi\|_{x,g}^*.$$
\end{defn}
\begin{rem}\label{comassf}
$\|\phi\|_{g}^*$ will be viewed as a pointwise function in this paper.
Generally it is merely continuous.
At a point $x$ where $\phi_x\neq0$,
\begin{equation*}
\begin{split}
\|\phi\|^*_{x,g}& =  \max \{\phi(\overrightarrow V_x): \overrightarrow V_x\ \text{is a simple $m$-vector  at } x \text{ with}\ 
\|\overrightarrow V_x\|_{g}=1\}\\
&=\max\{1/{\|\overrightarrow V_x\|_{g}}: \overrightarrow V_x\ \text{is\ a\ simple\ $m$-vector\ at } x \text{ with}\ \phi(\overrightarrow V_x)=1\}\ \\
&=
1/\min\{\|\overrightarrow V_x\|_{g}:\overrightarrow V_x\ \text{is\ a\ simple\ $m$-vector\ at } x \text{ with}\ \phi(\overrightarrow V_x)=1\}.
\end{split}
\end{equation*}
\end{rem}

\begin{defn}
Let  $(\mathscr E'_*(X),d)$ denote the dual complex of the {\em de Rham} complex %$\mathscr E^*(X)$
of $X$.
Elements of $\mathscr E'_k(X)$ are $k$-dimensional {\em de Rham} $\mathrm{\bold{currents}}$ (with compact support)
and $d$ is the adjoint of exterior differentiation.
%$dT(\phi)=T(d\phi)$
\end{defn}

%\begin{defn}
%A submanifold is called $\mathrm{\mathbf{calibratable}}$, if there exists some calibration form calibrating it (with respect to the given metric).
%\end{defn}

%Suppose $\phi$ calibrates a compact submanifold $M$ in $(X,\phi,g)$.
%If $M'$ is another compact oriented submanifold in the same $\mathbb R$-homology class of $X$,
%\begin{equation*}
%\text{Vol}_g(M)=\int_M\phi=\int_{M'}\phi\leq\ \text{Vol}_g(M'),
%\end{equation*}
%i.e., $M$ is mass-minimizing among ``smooth" objects in its homology class.
%%(The middle equality is still true for general case by {Theorem \ref{hl}}.)
%Actually this can be naturally
%generalized to the topological dual space $\mathfrak D_c$ of the space $\mathcal{E}$ of smooth forms, whose elements are called {\it de Rham} \textbf{currents} (with compact support).

%In fact there are two natural dualities $\mathfrak D_c\leftrightarrow\mathcal{E}\  \mathrm{and} \ \mathfrak D\leftrightarrow \mathcal{E}_c.$
%To avoid additional requirements on calibrations
% we use the first in calibrated geometry.

\begin{defn}\label{dmass}
Let $T$ be a {\em de Rham} $m$-current in $(X,g)$.
The ${\bold{mass}}$ of T is
$
\mathrm{\mathbf{M}}(T) = \sup \{T(\psi):\psi\ smooth\ m\text{-}form\ with\  \sup_X\|\psi\|_{g}^*\leq 1\}. 
$
\end{defn}
When $\mathrm{\mathbf{M}}(T)<\infty$, $T$ determines a unique {\it Radon} measure $\|T\|$ characterized by
$ \int_X f\cdot d\|T\|=\sup\{T(\psi): \|\psi\|_{x,g}^*\leq f(x) \}$
for any nonnegative continuous function $f$ on $X$.
Therefore $\mathrm{\mathbf{M}}(T)=\|T\|(X)$.
Moreover, the {\it Radon}-{\it Nikodym} Theorem asserts the existence of a $\|T\|$ measurable tangent $m$-vector field $\overrightarrow T$ a.e. with 
vectors $\overrightarrow T_x \in \Lambda^m T_xX$ of unit length in the dual norm of the comass norm,
% in {definition \ref{comassdef}},
satisfying
\begin{equation}\label{current}
T(\psi)= \int_X\psi_x(\overrightarrow {T_x})\ d \|T\|(x)\ \  \text{for any smooth $m$-form }\psi,
\end{equation}
\text{or\ briefly}
$T = \overrightarrow T\cdot \|T\|\ a.e.\ \|T\|.$
%Consequently two specific representatives of $T$ coincide except on a $\|T\|$ measure $0$ set. 

When $T$ enjoys local finite mass,
one can get {\it Radon} measure $\|T\|$ and decomposition (\ref{current}) as well.
\begin{defn}
Let $\mathbf{spt}(f)$ be the support of $f$ where $f$ is a function.
For a current $T$, let $U_T$ stand for the largest open set with $\|T\|(U_T)=0$.
Then the support of $T$ is denoted by $\mathbf{spt}(T)=U_T^c$.
\end{defn}
\begin{defn}
Let $\mathbb M_k(X)=\{T\in\mathscr E'_k(X): \mathrm{\mathbf{M}}(T)<\infty\}$.
Then $N_k(X)=\{T\in\mathbb M_k(X): dT\in\mathbb M_{k-1}(X)\}$ is the space of $k$-dimensional $\mathrm{\mathbf{normal}}$ currents.
\end{defn}
      \begin{rem}
      We view a current in $\mathbb M_k$ as a functional over smooth $k$-forms
      not a specific representative of generalized distribution.
      \end{rem}
Note that $(N_*(X),d)$ form a complex.
Recalling the natural isomorphisms established by {de Rham} and Federer and Fleming:
$$H_*(\mathscr E'_*(X))\cong H_*(X;\mathbb R)\cong H_*(N_*(X))$$
we identify these three homology groups the same.

\begin{defn}\label{calibration}
A smooth form $\phi$ on $(X,g)$ is called a $\mathrm{\mathbf{calibration}}$ if 
$
\sup_{X}\|\phi\|_{g}^*= 1$
and
$
d\phi=0.
$
The triple $(X,\phi,g)$ is called a $\bold{calibrated}$ $\bold{manifold}$. 
If $M$ is an oriented (singluar) submanifold
with $\phi|_{M}$ equal to the volume form of $M$,
then
$(\phi,g)$ is a $\mathrm{\mathbf{calibration\ pair}}$ of $M$ on $X$.
We say $\phi$
$\mathrm{\mathbf{calibrates}}$ $M$
and $M$ $\mathrm{\mathbf{can\ be\ calibrated}}$ in $(X,g)$.
\end{defn}

\begin{defn}\label{calibratable}
Let $\phi$ be a calibration on $(X,g)$.
We say that a current $T$ of local finite mass is $\bold{calibrated}$ by $\phi$, if 
$\phi_x(\overrightarrow T_x)=1\ a.a.\ x\in X\ \text{for}\ \|T\|.$
\end{defn}

The following is the fundamental theorem of calibrated geometry in \cite{HL2}.

\begin{thm}[Harvey and Lawson]\label{hl}
If $T$ is a calibrated current with compact support in $(X,\phi,g)$ and $T'$ is any compactly supported current homologous to $T$(i.e., $T-T'$ is a boundary and in particular $dT=dT'$), then
\begin{equation*}
\mathrm{\mathbf{M}}(T)\leq  \mathrm{\mathbf{M}}(T')
\end{equation*}
with equality if and only if $T'$ is calibrated as well.
\end{thm}

\begin{rem}\label{coneccc}
Let $M$ be an oriented compact smooth submanifold.
Then by continuity method, the current $[[M]]=\int_M\cdot\ $ is calibrated if and only if
$M$ is calibrated.
\end{rem}

We shall use certain properties of comass. 
Especially, {Lemma \ref{fperp}} is crucial to our methods and {Lemma \ref{CCGP}} gives a concrete control of comass for metric gluing.
\begin{lem}\label{norm}
For any metric $g$, $m$-form $\phi$ and positive function $f$ on $X$,
$$
\|\phi\|_{f\cdot g}^*=f^{-\frac{m}{2}}\cdot \|\phi\|_g^*.
$$
\end{lem}
\begin{prooff}
By the formula in {Remark \ref{comassf}}.
\end{prooff}
\begin{lem}\label{big}
For any $m$-form $\phi$ and metrics $g'\geq g$ on $X$, we have
$$
\|\phi\|_{g'}^*\leq\|\phi\|_g^*.
$$ 
\end{lem}
\begin{prooff}
By the definition of comass.
\end{prooff}
 \begin{lem}[Comass Control for Gluing Procedure]\label{CCGP}
For any $m$-form $\phi$, positive functions $a$ and $b$, and metrics $g_1$ and $ g_2$,
it follows
\begin{equation}\label{gluemetrics}
 \|\phi\|^*_{ag_1+bg_2}\leq\frac{1}{\sqrt{a^m\cdot \frac{1}{\|\phi\|^{*2}_{g_1}}+b^m\cdot \frac{1}{\|\phi\|^{*2}_{g_2}}}}
 \end{equation}
 where $\frac{1}{0}$ and $\frac{1}{+\infty}$ are identified with $+\infty$ and $0$ respectively.
 \end{lem}
 \begin{prooff}
The statement is trivial where $\phi=0$.
Now consider a point $x$ at where $\phi_x\neq0$.
In the subspace spanned by a simple $m$-vector $\overrightarrow V_x$,
there exists an orthonormal basis $(e_1,\ \cdots, e_m)$ of $g_1$, under
which $g_2$ is diagonalized as $diag(\lambda_1,\ \cdots ,\lambda_m)$
for some $\lambda_i>0$.
Suppose $\overrightarrow V_x=te_1\wedge\cdots\wedge e_m$, then
\begin{equation}\label{inq}
\begin{split}
\|\overrightarrow V_x\|_{ag_1+bg_2}^2 &= t^2(a+b\lambda_1)\cdots(a+b\lambda_m)\\
&=\  \ t^2[a^m+\cdots+b^m{\bold \Pi}\lambda_i]\\
&\geq \ \ \quad t^2a^m+t^2b^m{\bold\Pi}\lambda_i\\
%&= a^m\|\overrightarrow V_x\|_{h_1}^2+b^m\| \overrightarrow V_x\|_{h_2}^2\\
&= a^m\|\overrightarrow V_x\|_{g_1}^2+b^m\|\overrightarrow V_x\|_{g_2}^2. 
\end{split}
\end{equation}
By {Remark \ref{comassf}},
%it is not hard to see that 
(\ref{inq}) implies \eqref{gluemetrics}.
\end{prooff}
\begin{lem}[Comass One Lemma]\label{fperp}
Suppose $(E, \pi)$ is a disc bundle over $M$ (as the zero section) and $g$ is a Riemannian metric defined on $E$.
Then each fiber is perpendicular to $M$ if and only if $\pi^*vol_{g_M}$ has comass one pointwise along $M$
where $vol_{g_M}$ means the volume form of $M$ induced by $g$. 
\end{lem}

\begin{prooff}
Fix a point $x$ on $M$. Take an oriented orthonormal basis $\{e_1,\ \cdots,\ e_m\}$ of $T_xM$. 
Then we have unique decompositions
$
e_i=\sin\theta_i\cdot a_i+ \cos\theta_i\cdot b_i
$
where $a_i$ is a unit vector perpendicular to $F_x$ $-$ the subspace of
the fiber directions in $T_xE$, $b_i$ is some unit vector in $F_x$,
and $\theta_i$ is the angle between $e_i$ and $F_x$ for $i=1,\ \cdots, m$.
Denote $vol_{g_M}$ by $\omega$. By the choice of $\{e_i\}$,
\begin{equation}\label{q1}
\begin{split}
{\ } \ \ \ \ \ \ \ \ 1&=\omega(e_1\wedge e_2\cdots\wedge e_m)\\
&=\pi^*\omega(e_1\wedge e_2\cdots\wedge e_m)\\
&=\pi^*\omega(\sin\theta_1\cdot a_1\wedge\cdots\wedge\sin\theta_m\cdot a_m)\\
&={\bold \Pi} \sin\theta_i\cdot \pi^*\omega(a_1\wedge a_2\cdots\wedge a_m).
\end{split}
\end{equation}
The third equality is due to the fact that elements of $F_x$ annihilate $\pi^*\omega$.
Since $\{a_i\}$ are of unit length,
$\|\pi^*\omega\|^*_{x,g}\geq 1, \ \forall x\in M.$
Combined with {Remark \ref{comassf}}, the equality holds if and only if $F_x \perp T_xM$.
\end{prooff}

Since the pullback of volume form is simple and well behaves, we want to extend its certain multiple to a global potential calibration form. In order to proceed,
we need global forms.
{\ }
\\
%%%%%subsection
\subsection{For Global Forms}\label{basic}
In the singular homology theory
the {\em Kronecker} product $<\cdot ,\cdot >$ between cochains and chains induces a homomorphism
\begin{equation*}\label{kappa}
\kappa\ : H^q(X; G)\rightarrow \ \text{Hom}_{\mathbb Z}(H_q(X;\mathbb{Z}),\ G) \text{ given by}
\end{equation*}
\begin{equation*}
\kappa\ ([z^q])([z_q])\triangleq <[z^q],[z_q]>
\end{equation*}
where $G$ is any Abellian group.
A classical result asserts that
$\kappa$ is surjective. 
When $G=\mathbb{R}$, by {\it de\ Rham} {Theorem},
$
\kappa\ : H_{dR}^q(X)\twoheadrightarrow \  \text{Hom}_{\mathbb{R}}(H_q(X;\mathbb{R}),\mathbb{R}).
$

Suppose $\{M_\alpha\}$ are mutually disjoint $m$-dimensional oriented connected compact submanifolds with
the represented homology classes $\{[M_\alpha]\}$ lying in one common side of some hyperplane through the zero of $H_m(X;\mathbb{R})$.
Then there exists a homomorphism 
$\digamma\in \text{Hom}_{\mathbb{R}}(H_m(X;\mathbb{R}),\mathbb{R})$
forwarding $\{[M_\alpha]\}$ to positive numbers.
Consequently, we have the following existence result.

\begin{lem}\label{formmore}
Suppose $X$ is a manifold
and $\{M_\alpha\}$ satisfy the above condition.
Then there exists a closed $m$-form $\phi$ on $X$ with
$
\int_{M_\alpha} \phi>0\ \text{ for  every } M_\alpha.
$
\end{lem}

\begin{rem}
When $\#\{M_\alpha\}=r<\infty$, the requirement becomes the
$\bold{convex\ hull\ condition}$ that
$\{\sum_{i=1}^rt_i[M_i]: \sum_{i=1}^rt_i=1\ \mathrm{and}\ t_i\geq 0\}$
in $H_m(X;\mathbb{R})$
does not contain the zero class.
\end{rem}
%\begin{rem}
%Since our analysis above is also valid for the case of countable components (to occur in $\S \ref{Xnonc}$ and $\S\ref{ccm}$),
%the same statement holds with no essential differences for that case as well.
%\end{rem}
{\ }
\\
\subsection{Bundle Structure around Submanifolds}\label{BSNS}
Given an oriented compact submanifold $M$ in $(X,g)$, consider its $\epsilon$-neighborhood $U_\epsilon%\triangleq \{x\in X: \text{dist}_g(x,M)\leq\epsilon\}
$.
When $\epsilon$ is small enough, 
the metric induces a disc-fibered bundle structure of $U_\epsilon$,
whose fiber is given by the exponential map restricted to normal directions along $M$.
%It is clear that the fibers are foliated by the distance function from $M$.
We call the orthogonal complement
to fiber directions in $TU_\epsilon$ \textbf{horizontal directions} and a
\textbf{horizontal change} means a smoothly varying modification on $g$
along horizontal directions.

By a strong deformation retraction from $U_{\epsilon}$ to $M$,
$
H^m(U_\epsilon;\mathbb{R})\cong H^m(M;\mathbb{R}).
$
Therefore for any $[\phi_1]$ and $ [\phi_2]\in$ $H^m(U_\epsilon;\mathbb{R})$
\begin{equation}\label{equt}
[\phi_1]=[\phi_2]\ \Leftrightarrow \int_M \phi_1=\int_M \phi_2.
\end{equation}

For an oriented properly embedded non-compact complete submanifold (without boundary), 
similar bundle structure occurs for
some suitable smooth positive function $\epsilon$ defined along the submanifold.
\\

Besides a manifold, a calibrated manifold consists of a balanced pair of a metric and a closed form. 
We shall do gluing procedures in both slots.
{\ }\\
\subsection{Gluing of Forms}\label{forms}
Suppose $M$ is an oriented compact connected $m$-dimensional submanifold
with $[M]\neq[0]\in H_m(X;\mathbb R)$.
By {\S \ref{basic}}
there exists a closed $m$-form $\phi$ on $X$ with 
$s=\int_M \phi> 0.$

Let $\pi_g$ be the projection map of the disc bundle in
{\S \ref{BSNS}}.
Then in $U_\epsilon$ 
\begin{equation}\label{ef}
\int_M \frac{s\cdot\pi^*_g\omega_M}{\text{Vol}_g(M)}=s=\int_M  \phi
\end{equation}
where
$\omega_M$ is the volume form of $M$.
Denote the integrand of the left hand side of \eqref{ef} by $\omega^*$.
By (\ref{equt}) $[\omega^*]=[\phi]$ in $H^m(U_\epsilon;\mathbb{R})$
which indicates
\begin{equation}\label{abouttoglue}
\phi=\omega^*+d\psi
\end{equation}
for some smooth $(m-1)$-form $\psi$ defined on $U_\epsilon$.
Now take $\Phi=\omega^*+d((1-\rho(\textbf{d}))\psi)$
where \textbf{d} is the distance function to $M$ and $\rho$ is given in the picture.
%********************************************* graph1************
 \begin{figure}[ht]
\centering
\includegraphics[scale=0.17]{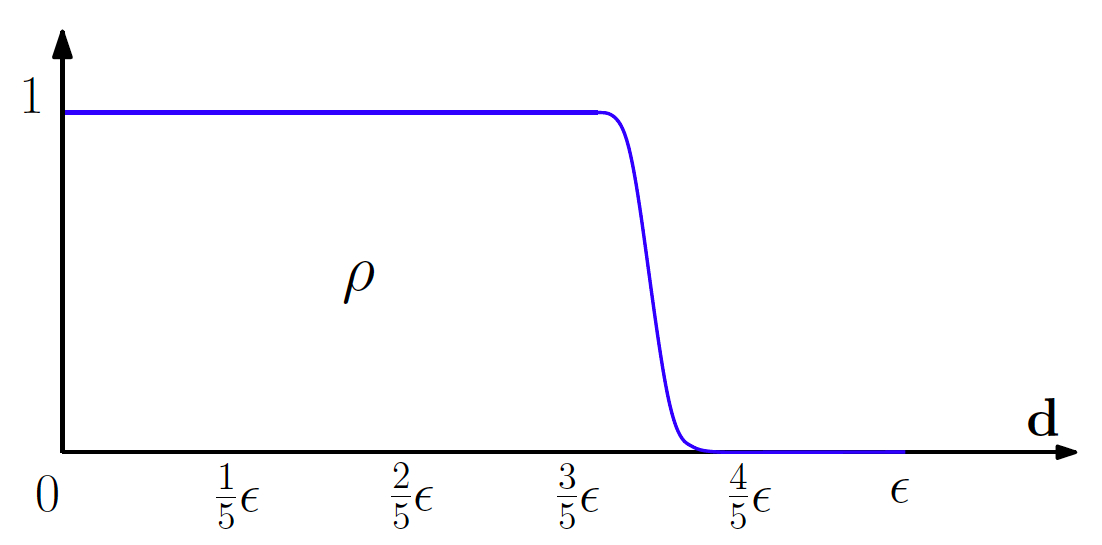}
\end{figure}

Clearly $\Phi$ extends to a closed smooth form on $X$ as follows.
\begin{equation*}
\Phi=
\begin{cases}\label{glueforms}
\omega^* & 0\leq  \textbf{d}\leq \frac{3}{5}\epsilon\\
\omega^* + d(1-\rho( \textbf{d}))\psi\ \ \ \ \ \ \ & \frac{3}{5}\epsilon <  \textbf{d} \leq \frac{4}{5}\epsilon\\
\phi & \frac{4}{5}\epsilon< \textbf{d}\\
\end{cases}
\end{equation*}

By {Lemma \ref{fperp}},
$\|\Phi\|_g^*=\|\omega^*\|^*_{ g} =\frac{s}{\text{Vol}_g(M)}$ along $M$.
\\
{\ }
\\
{\ }
\\
{\ }

% Section 3

\section{Horizontal Change of Metric}\label{H} % Main section title

\label{Section3} % For referencing the section elsewhere, use \ref{Section1} 

\subsection{Gluing of Metrics}\label{metrics}
In this subsection we glue metrics for the horizontal change.
By the bundle structure of $U_\epsilon$, we define
\begin{equation*}\label{hcm}
\bar g=(\frac{s}{\text{Vol}_g(M)})^{\frac{2}{m}}\pi_g^*(g_M)\oplus 
g_F
\end{equation*}
where $
g_F$ is any smooth metric along fiber directions and 
$\pi_g^* (g_M)(v_q,v_q')=  g_M(\pi^g_*(v_q), \pi^g_*(v_q'))$
where $\pi^g_*$ is the push-forwarding map of $\pi_g$.
\\

\begin{prop}\label{property}
$\|\omega^*\|^*_{\bar g} =1$ on 
$U_\epsilon$.
\end{prop}
\begin{prooff}
Fix a point $q$ in $U_\epsilon$.
By {Remark \ref{comassf}} one only needs to verify that
\begin{align}\label{min}
\min\{\|W\|_{\bar g}:W\ \text{is a simple } m\text{-vector at } q \text{ with }  
\omega^*(W)=1\}
\end{align}
equals one.
Suppose $\overline W$ realizes the minimum of \eqref{min} and
 it decomposes into a purely horizontal part and a fiber-involving part:
$\overline W=\overline W^h+ \overline W^\nu.$
By Lemma \ref{norm} and definitions of $\omega^*$ and $\bar g$,
\begin{equation*}
\begin{split}
1&=\omega^*(\overline W)
=\frac{s\pi^*_g \omega_M(\overline W)}{\text{Vol}_g(M)}
=\frac{s\pi^*_g \omega_M(\overline W^h)}{\text{Vol}_g(M)}
=\frac{s\omega(\pi^g_*\overline W^h)}{\text{Vol}_g(M)}\\
&=\frac{s\|\pi^g_*\overline W^h\|_{g_M}}{\text{Vol}_g(M)}
=\|\pi^g_*\overline W^h\|_{(\frac{s}{\text{Vol}_g(M)})^{\frac{2}{m}}g_M}
=\|\overline W^h\|_{\bar g}\leq \|\overline W\|_{\bar g}.
\end{split}
\end{equation*}
The equality holds if and only if $\overline W$ is purely horizontal.
Hence \eqref{min} is indeed one.
\end{prooff}
\begin{rem}\label{uniquevec}
Due to the dimension reason, $\overline W$ is the unique oriented horizontal $m$-vector with $\omega^*(\overline W)=1$.
In particular
$\|\omega^*\|^*=\frac{1}{\|\overline W\|}$ for any metric.
\end{rem}

We construct a new metric $\tilde g= \tilde g^h\oplus\tilde g^{\nu}$
based on $g=g^h\oplus g^\nu$ by gluing metrics  along horizontal and fiber directions respectively:
\begin{align}
 \tilde g ^h&=\sigma^{\frac{1}{m}}
 ((\frac{s}{\text{Vol}_g(M)})^{\frac{2}{m}}+\textbf{d}^2)
 \pi_g^* (g_M)+(1-\sigma)^{\frac{1}{m}}
 \alpha g^h, \text{ and} \label{fhcm}\\
 \tilde g ^v&=\sigma g_F + (1-\sigma) \alpha g^v,\label{2fhcm}
\end{align}
where \textbf{d} is as previous,
$\alpha$ is a smooth positive function (to be determined below) on $X$,
and $\sigma=\sigma$(\textbf{d}) is given in the picture.
%********************************************* graph2**************
 \begin{figure}[ht]
\begin{center}
\includegraphics[scale=0.22]{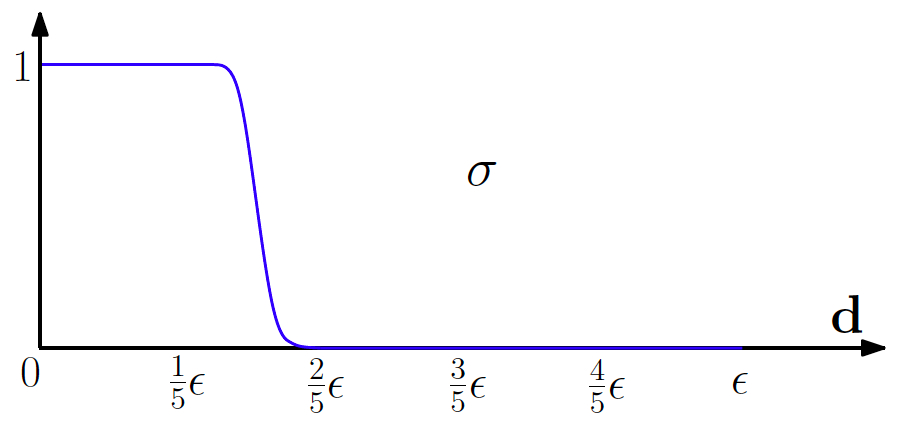}
\end{center}
%\caption{Graph of $\sigma$.}
\end{figure}

Let $g_F$ be $\alpha g^v$.
Then
$
 \tilde g ^v=\alpha g^v.
$
Altough expressions \eqref{fhcm} and \eqref{2fhcm} are valid in $U_\epsilon$ only, $\tilde g$ can be extended to a global metric by defining $\tilde g=\alpha g$ on $(U_\epsilon)^c$.
Choose an appropriate $\alpha$ such that
\begin{equation}\label{globalcontrolofalpha}
\|\Phi\|^*_{\alpha g}<1\text{ on } X,
\end{equation}
which implies
\begin{equation}\label{localcontrolofalpha}
\alpha^{{m}}\|\overline W\|_g^2>1 \text{ on } U_{\frac{3}{5}\epsilon}.
\end{equation}
Then
{\small
$$\|\Phi\|^*_{\tilde g}=
\begin{cases}
\ \ \ \ 1& \ 0= \textbf{d}\\
\ \|\omega^*\|^*_{\tilde g}\ \ <\ \ \|\omega^*\|^*_{\bar g}=1
%\ \ =\ \frac{1}{\|\overline W\|_{((\frac{s}{\text{Vol}_g(M)})^{\frac{2}{m}}+\textbf{d}^2){\pi_g^* (g_M)}}}\ \ <\ 1
& \ 0< \textbf{d} \leq \frac{1}{5}\epsilon\\
\\
\begin{array}{lll}
\|\omega^*\|^*_{\tilde g}
&=&\|\omega^*\|^*_{
%[\sigma^{\frac{1}{m}}((\frac{s}{\text{Vol}_g(M)})^{\frac{2}{m}}+\textbf{d}^2){\pi_g^* (g_M)}+\alpha(1-\sigma)^{\frac{1}{m}}g^h]
\tilde g^h
\oplus \tilde g^v}\\
&\leq& \frac{1}{\sqrt{\sigma\|\overline W\|^2_{((\frac{s}{\text{Vol}_g(M)})^{\frac{2}{m}}
+\textbf{d}^2){\pi_g^* (g_M)}}
+\alpha^m(1-\sigma)\|\overline W\|^2_{g^h}}}\\
&<&1\
\end{array}
&\frac{1}{5}\epsilon \leq  \textbf{d}\leq \frac{2}{5}\epsilon
\\
\ \|\Phi\|^*_{\alpha g}\ \ <\ \ 1\ &   \frac{2}{5}\epsilon \leq  \textbf{d}%\leq \epsilon
\end{cases}$$
}\\
The first inequality is by Proposition \ref{property},
the second from {Lemma \ref{CCGP}} and Remark \ref{uniquevec},
the third by Proposition \ref{property} and (\ref{localcontrolofalpha}), %when $\frac{1}{5}\epsilon \leq d\leq \frac{2}{5}\epsilon$
and the last due to (\ref{globalcontrolofalpha}). %when $\frac{2}{5}\epsilon \leq d$.
In summary, we obtain a calibration pair $(\Phi, \tilde g)$
with
$\bold{spt}(\|\Phi\|^*_{\tilde g}-1)=M$.

\begin{rem}\label{alpha}
If $X$ is compact, then $\alpha$ can be chosen as a sufficiently large constant 
to guarantee (\ref{globalcontrolofalpha}).
It is clear that the same constant $\alpha$ still works when one shrinks the neighborhood $U_\epsilon$.
\end{rem}
%(Maybe some convergent status for metrics could be made!!!!!!)
%%%%%%%%%%%%%%%%%%%%%%%%%%%%%%%%%%%%%%%%%%%%%%%
{\ }\\
\subsection{Compact $X$ and Connected $M$}
\begin{defn}\label{strongC}
A submanifold $M$ is $\mathrm{\mathbf{strongly}}$ $\mathrm{\mathbf{calibrated}}$ (or $\mathrm{\mathbf{tamed}}$) in a calibrated manifold $(X,\phi,g)$, 
if each connected component of $M$ is calibrated by $\phi$ (either $\phi$ or $-\phi$) and $\bold{spt}(\|\Phi\|^*_{g}-1)=M$. 
\end{defn}
\begin{thm}\label{a}
Suppose $(X,g)$ is a {\tt compact}
%%%%%%%%%%%%%%%oriented closed compact 
Riemannian manifold and $M$ 
is an oriented compact connected $m$-dimensional submanifold
representing a nonzero class of $H_m(X;\mathbb R)$.
Then for any open neighborhood $U$ of $M$,
a new metric $\hat g$ can be constructed by a horizontal change of $g$
supported in $U$
such that $M$ is strongly calibrated by some calibration $\hat \Phi$ in $(X, \hat g)$.
\end{thm}
\begin{prooff} 
By the compactness of $M$, there exists a small positive number $\epsilon$ to guarantee 
the disc bundle structure of $U_\epsilon$ in {\S \ref{BSNS}} and $U_\epsilon\subseteq U$.
By {\S \ref{metrics}} and Remark\ \ref{alpha},
a new metric $\tilde g$ can be constructed by a horizontal change in $U_\epsilon$ with 
$\tilde g=\alpha g$ on $X-U_\epsilon$, where $\alpha$ is a large constant satisfying (\ref{globalcontrolofalpha}).
Define $\hat g\triangleq  \alpha^{-1}\tilde g$ and $\hat \Phi \triangleq  \alpha ^{-\frac{m}{2}}\Phi$.
By {Lemma \ref{norm}}, 
$M$ is strongly calibrated in $(X, \hat \Phi, \hat g)$.
Moreover $\hat g$ equals $g$ on $X-U_\epsilon$.
\end{prooff}
%\begin{rem}\label{consistant}
%Once $\phi$ is chosen, our construction shows
%$$\mathrm{Vol}_{\tilde g}(M)=\int_M\phi$$
%no matter which initial metric we start with.
%\end{rem}

In fact metrics produced in the proof make $M$ more than being minimal.

\begin{prop} 
$M$ is totally geodesic under the constructed metric $\hat g$.
\end{prop}
\begin{prooff}
Note that the projection map is length-shrinking.
\end{prooff}
{\ }\\

\subsection{Compact $X$ and Non-Connected $M$}\label{Mnonc}
\begin{defn}\label{collection}
A family $\mathfrak M$ of mutually {disjoint} connected oriented compact submanifolds of a manifold $X$ (not necessarily compact) is called a $\mathrm{\mathbf{mutually}}$ $\mathrm{\mathbf{disjoint}}$ $\mathrm{\mathbf{collection}}$
and an element of $\mathfrak M$ is a 
$\mathrm{\mathbf{component}}$.
The (nonempty) subset $\mathfrak M_k$ of all components of dimension $k$ is its $\mathrm{\mathbf{k\text-level}}$.
When $\mathfrak M$ consists of finitely (or countably) many components,
we call it a $\mathrm{\mathbf{finite}}$ (or $\mathrm{\mathbf{countable}}$) 
$\mathrm{\mathbf{collection}}$.
If $\mathfrak M$ has only one level in dimension $m$,
then it is called an $\mathrm{\mathbf{m\text-collection}}$.
Let $\underline{\mathfrak M_k}$ denote the union submanifold of the components of $\mathfrak M_k$ and define $\underline{\mathfrak M}\triangleq\bigcup_k\underline{\mathfrak M_k}$ as a set.
\end{defn}

In this subsection we assume that $X$ is \underline{{compact}}.
%********************************************* Theorem 3:************************
\begin{thm}\label{mainthm1}
Let $\mathfrak M=\{M_i\}_{i=1}^s$ be a finite mutually disjoint $m$-collection satisfying the convex hull condition.
Then for any open neighborhood $U$ of $\underline{\mathfrak M}$,
a new metric $\hat g$ can be constructed by a horizontal change of $g$ supported in $U$
such that there exists a calibration $m$-form $\hat\Phi$ on $(X,\hat g)$
and every nonzero current $T= \sum_{i=1}^st_i[[M_i]]$ with $t_i\geq0$
is calibrated in $(X,\hat \Phi, \hat g)$.
Consequently $T$ is mass-minimizing in $[T]$ with $\mathrm{\mathbf{M}}(T)=\sum_{i=1}^st_i \mathrm{Vol}_{\hat g}(M_i)$.
\end{thm}

If we just require that
each component represents a nonzero class,
then there exists some hyperplane $\mathcal{P}_m$ through zero in $H_m(X;\mathbb{R})$
avoiding all classes $\{[M_i]\}_1^s$. $\mathcal{P}_m$ divides the space into two open chambers.
By reversing orientations of components of $\mathfrak M_m$ in one chamber,
we get a new collection satisfying the convex hull condition.

%Since mass is independent of the choice of orientations, the following is an immediate consequence.
\begin{cor}\label{maincor1}
%********************************************* Theorem 2*********\\
Suppose $\mathfrak M$ is a finite mutually disjoint $m$-collection in $(X, g)$
and each component represents a nonzero class in $H_m(X;\mathbb{R})$.
Then for any open neighborhood $U$ of $\underline {\mathfrak M}$,
 a new metric $\hat g$ can be constructed by a horizontal change of $g$ supported in $U$
such that $\underline {\mathfrak M_m}$ can be tamed in $(X, \hat g)$.
\end{cor}
\begin{pf}
By assumption there exists
a smooth closed $m$-form $\phi$ on $X$
with $\int_{M_i}\phi>0$ for every $M_i$. 
Since $\mathfrak M_m$ is a finite collection,
one can choose $\epsilon>0$
such that %(as described in {\S \ref{BSNS}})
$\{U_\epsilon(M_i)\}$ are mutually
disjoint with bundle structures in {\S \ref{BSNS}} and contained in $U$.
Gluing $\phi$ with local forms as in \S\ref{forms} on $\{U_\epsilon(M_i)\}$,
we get $\Phi$.
Using a large constant $\alpha$ satisfying \eqref{globalcontrolofalpha} for the metric gluing,
one can obtain a new metric $g_\alpha$ under which
$\Phi$ becomes a calibration.
Let $\hat g\triangleq  \alpha^{-1}g_{\alpha}$ and $\hat \Phi \triangleq  \alpha ^{-\frac{m}{2}}\Phi$.
Then every $M_i$ is calibrated by $\hat \Phi$ in $(X,\hat g)$.

Since $\{M_i\}$ are submanifolds, the measure $\|[[M_i]]\|$ is 
$\mathcal{H}|_{M_i}$, the {\it Hausdorff} $m$-measure restricted
to $M_i$ in $(X, \hat g)$.
When $T=\sum_{i=1}^st_i[[M_i]]$ with $t_i\geq 0$, 
we have $\|T\|=\sum_i^st_i 
\mathcal{H}|_{M_i}$ and therefore
$
\textbf{M}(T)%&
=
\int_X\|\overrightarrow T\|_{\hat g}\ d\|T\|
%=
%\sum t_i \int_{M_i}\|\overrightarrow M_i\|_{\hat g}\ d\mathcal{H}|_{M_i}\\&
=\sum t_i \int_{M_i}d\ vol_{\hat g}=\sum t_i \text{Vol}_{\hat g}(M_i).
$
\end{pf}
\begin{rem}\label{yiyang}
If
$\sum s_i [M_i]=\sum t_j [M_j] \ for\ some\ s_i,t_j\geq 0,$
then we have
$\sum s_i\mathrm{Vol}_{\hat g}(M_i)=\sum t_j\mathrm{Vol}_{\hat g}(M_j)$
automatically.
\end{rem}

{\ }\\
\subsection{Compact $X$ and $\mathfrak M$ with Several Levels}\label{Mdiffd}
Assume $\mathfrak M$ is a finite mutually disjoint collection, not necessarily of a single level, in a \underline{{compact}} Riemannian manifold $(X,g)$.
Results similar to {Theorem \ref{mainthm1}} and {Corollary \ref{maincor1}} can be gained.
%Since the methods are almost the same, here we only provide a proof for {\it Corollary \ref{dimdiffcor}} below.
\begin{thm}\label{dimdiffthm}
%********************************************* Theorem 3:************************

Suppose $\mathfrak M$ is a finite mutually disjoint collection
and each $k$-level $\mathfrak M_k$ satisfy the convex hull condition.
Then for any open neighborhood $U$ of $\underline{\mathfrak M}$,
a new metric $\hat g$ can be constructed by a horizontal change of $g$ supported in $U$
such that there exist
a family of calibrations 
$\{\hat \Phi_k\}$ in $(X,\hat g)$
and every nonzero current $T=\sum_{i=1}^{s_k} t_i[[M_i]]$
with $M_i\in\mathfrak M_k$ and $t_i\geq0$
is calibrated by $\hat \Phi_k$. 
Consequently $T$ is mass-minimizing
in $[T]$
with $\mathrm{\mathbf{M}}(T)=\sum_{i=1}^st_i \mathrm{Vol}_{\hat g}(M_i)$.
\end{thm}

\begin{cor}\label{dimdiffcor}
%******************
Let $\mathfrak M$ be a finite mutually disjoint collection in $(X, g)$
with each component representing a nonzero class in the $\mathbb{R}$-homology of $X$.
Then for any open neighborhood $U$ of $\underline {\mathfrak M}$,
 a new metric $\hat g$ can be constructed by a horizontal change of $g$ supported in $U$
such that every $\underline{\mathfrak M_k}$ can be tamed in $(X, \hat g)$.
\end{cor}
\begin{pf2}
Without loss of generality,
assume $\mathfrak M=\{A^a,B^b\}$ with $a>b$. 
Take $\epsilon>0$ to guarantee bundle structures of $U_\epsilon(A)$ and $U_\epsilon(B)$ in {\S \ref{BSNS}},
$U_\epsilon(A)\cap U_\epsilon(B)=${ {$\O$}} and $U_\epsilon(A\bigcup B)\subseteq U$.
Suppose one gets $\Phi$ for $A$ as in \S\ref{forms}.
Then $\Phi=d\theta$ in $U_{\epsilon}(B)$. 
Set $\tilde \Phi=\Phi-d(\tilde\rho(\textbf{d}) \theta)$ where $\tilde\rho$ is given in the picture.
Then $\tilde 
\Phi$ is zero on $U_{\frac{4}{5}\epsilon}(B)$.
There exists a function $\alpha_a$
with value one on $\{\tilde \rho=1\}$
satisfying \eqref{globalcontrolofalpha}.
Based on $\tilde \Phi, g$ and $\alpha_a$, one can get a metric $\tilde g$ as in {\S \ref{metrics}}
such that $A$ is calibrated by $\tilde \Phi$ under $\tilde g$.
%*************************************************graph5**********************
 \begin{figure}[ht]
\begin{center}
\includegraphics[scale=0.2]{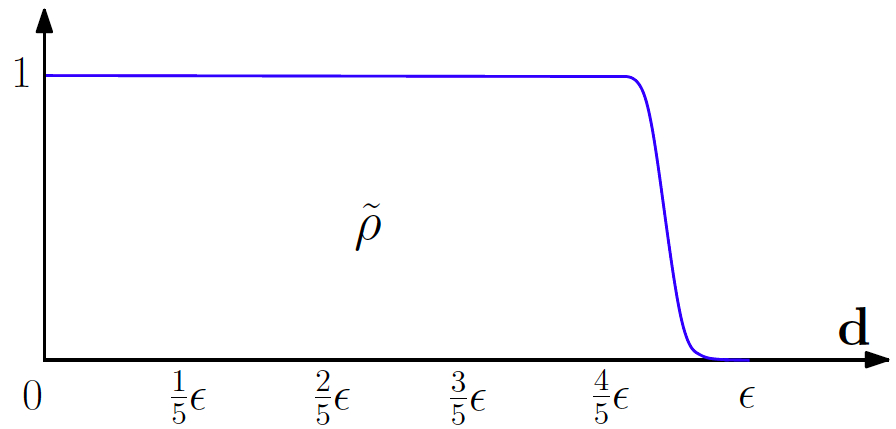}
\end{center}
%\caption{Graph for $\tilde \rho$.}
\end{figure}

By the compactness of $X$, one can choose a $b$-form $\psi$ with $\int_B\psi>0$ and $\|\psi\|_{\tilde g}^*<1$ on $X$.
Then by taking $\alpha_b\equiv1$, one can get $(\Psi,\hat g)$ as in \S\ref{metrics}.
Now $\tilde \Phi$ and $\Psi$
calibrate $A$ and $B$
respectively in $(X,\hat g)$.
\end{pf2}
%\begin{rem}
%Combined with the proof of {\it Theorem \ref{mainthm1}}, it is not hard to obtain the claims in {\it Theorem \ref{dimdiffthm}}.
%\end{rem}
\begin{rem}
Another proof is to choose a proper constant factor $\alpha$ for each level without changing
potential calibration forms.
However it does not work for the conformal change in {\S \ref{ccm}}.
\end{rem}

\begin{prop}\label{wps}
For the resulted metric $\hat g$,
$\mathrm{dist}_{\hat g}(M_i,M_j)= \mathrm{dist}_{g}(M_i,M_j)$.
\end{prop}%*************************************************graph4**********************
 \begin{figure}[ht]
\begin{center}
\includegraphics[scale=0.17]{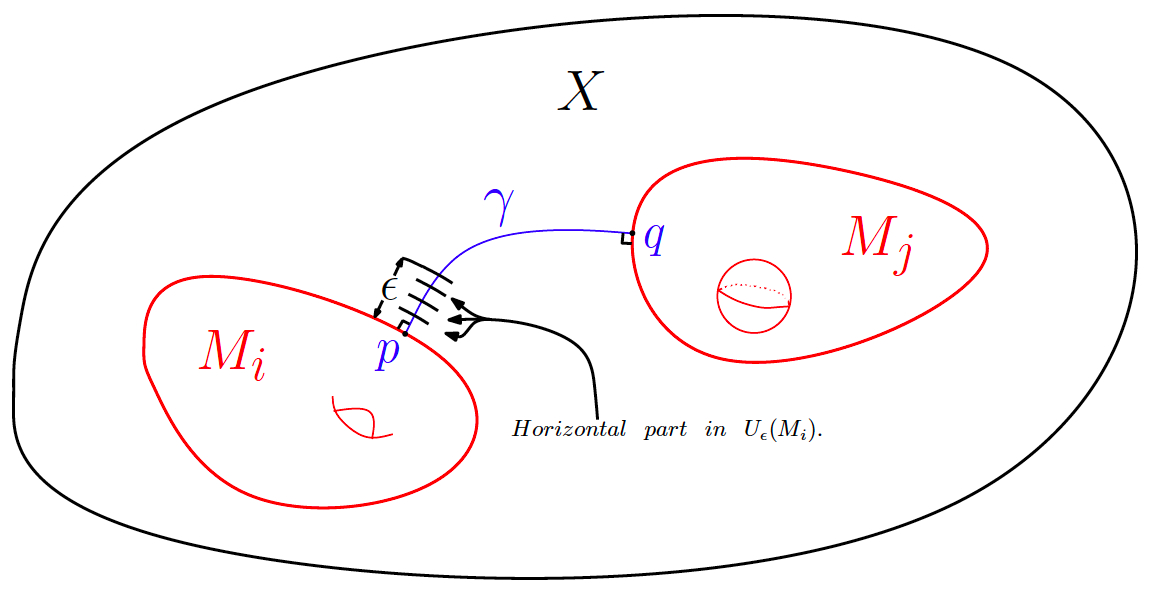}
\end{center}
%\caption{Graph for same distance.}
\end{figure}
\begin{prooff}
Suppose
$\gamma$ is a geodesic segment 
realizing $\mathrm{dist}_{g}(M_i,M_j)$.
According to {\S \ref{BSNS}}, $\gamma\cap U_\epsilon(M_i)$
is contained in some single fiber of $U_\epsilon(M_i)$ and similarly around $M_j$.
Since $\hat g$ and $g$ are the same along fiber directions, it follows that
$
\mathrm{dist}_{\hat g}(M_i, M_j)\leq l_{\hat g}(\gamma)=l_{g}(\gamma)=\mathrm{dist}_g(M_i, M_j)
$
where $l$ stands for the length functional.
{Appendix\ \ref{C}} shows that
the foliations induced by $g$ and $\hat g$ of $U_\epsilon(M_i)$ coincide.
Hence the same kind of argument leads to the opposite inequality.
\end{prooff}
{\ }\\

\subsection{Non-Compact $X$}\label{Xnonc}
We want to extend our gluing methods to the case of non-compact ambient manifolds.
However
in contrast to {Remark \ref{alpha}} for the compact case, 
generally a constant function $\alpha$ cannot meet the need of (\ref{globalcontrolofalpha}). 
For example, consider the surface $X$ obtained by rotating the graph of $y=e^x$ around the x-axis
with the induced metric from Euclidean $\mathbb{R}^3$. 
Take $\gamma$ as the circle
corresponding to $x=0$.
Obviously there is no way to
make $\gamma$ homologically minimal
through a local change of metric.
%With necessary modification, similar ideas apply to a collection in a non-compact Riemannian manifold $(X,g)$.
%The proof of {\it Theorem \ref{mm}} given below also works for {\it Theorem \ref{dimdiffthm}}.
\begin{thm}\label{mm}
Suppose $\mathfrak M$ is a finite mutually disjoint  collection
and each component represents a nonzero class in the $\mathbb{R}$-homology of $X$.
Then for any open neighborhood $V$ of $\underline {\mathfrak M}$ with a compact closure $\overline V$,
a new metric $\hat g$ can be constructed by a horizontal change of $g$ supported in $V$ and a conformal change in $X-V$
 such that every $\underline {\mathfrak M_k}$ can be tamed in $(X,\hat g)$.  
\end{thm}
\begin{prooff}
One can modify the proof of Theorem \ref{dimdiffthm}.
based on the compactness of $\overline V$.
\end{prooff}

In order to apply gluing techniques to a countable mutually disjoint collection, we introduce the following definition.
\begin{defn}\label{neat}
Let $\mathfrak M=\{M_i\}_{i=1,2,\cdots}$ be a countable mutually disjoint collection in $X$.
If the set $\bigcup_{i\neq j}M_i$ is closed for any positive integer $j$, 
then we call $\mathfrak M$ a $\mathrm{\mathbf{neat\ mutually\ disjoint}}$  $\mathrm{\mathbf{collection}}$.
\end{defn}

The neatness implies the existence of $\epsilon_i>0$ such that $U_{\epsilon_i}(M_i)$ are
mutually disjoint.
Hence we get the following.
%However the method fails to descend further from a level of infinitely many components.

\begin{thm}\label{Xnoncthm}
Suppose $\mathfrak M$ is a neat mutually disjoint collection
and each component represents a nonzero class in the $\mathbb{R}$-homology of $X$.
In addition, assume every level of $\mathfrak M$ consist of finite components except the lowest (nonempty) level.
 Then for any open neighborhood $U$ of $\underline {\mathfrak M}$,
 a new metric $\hat g$ can be constructed by a horizontal change of $g$ supported in $U$
 and a conformal change in the complement of the support of the horizontal change
such that every $\underline{\mathfrak M_k}$ can be tamed in $(X,\hat g)$.  
\end{thm}
{\ }
\\

% Section 4

\section{Conformal Change of Metric}\label{ccm}

\subsection{Parallel Results to Section \ref{H}}\label{pr}
An interesting kind of metric changes is the conformal change,
the simplest way to vary a metric.
Let us first glance at a basic case to understand our gluing method for conformal change. 
\begin{thm}\label{mainthm2}
Suppose $(X,g)$ is a {\tt compact}
Riemannian manifold and
$M$ 
is an oriented compact connected $m$-dimensional submanifold
representing a nonzero class in $H_m(X;\mathbb{R})$. Then for any open neighborhood $U$ of $M$,
a new metric $\hat g$ can be constructed by a conformal change of $g$ supported in $U$ such that
$M$ is strongly calibrated by some calibration $\hat\Phi$ in $(X, \hat g)$.
\end{thm}
\begin{prooff}
Since $\Phi$ in {\S \ref{basic}}
is pointwise a multiple of $\pi_g^*\omega$ in $U_{\frac{3}{5}\epsilon}(M)$,
$\|\Phi\|^*_g$ is smooth on $U_{\frac{3}{5}\epsilon}(M)$.
Let $g'= (\|\Phi\|^*_g )^{\frac{2}{m}}g$.
Then $\|\Phi\|^*_{g'}=1 \text{ on } U_{\frac{3}{5}\epsilon}(M).$
Take
\begin{equation}\label{mgofcc}
\tilde g=\sigma^{\frac{1}{m}}(1+ \textbf{d}^2)g'+\alpha(1-\sigma)^{\frac{1}{m}}g
\end{equation}
where $\sigma$ and $\textbf{d}$ are the same as in {\S \ref{metrics}}
and $\alpha$ is a large constant 
satisfying (\ref{globalcontrolofalpha}).
Set $\hat g\triangleq\alpha^{-1}\tilde g$ and $\hat \Phi \triangleq \alpha ^{-\frac{m}{2}}\Phi$.
Then $M$ is strongly calibrated by $\hat\Phi$ under $\hat g$.
Furthermore,
$\hat g$ is conformal to $g$
and $\hat g=g$ on $X-U_\epsilon(M)$.
\end{prooff}

The same local gluing ideas and elimination tricks on calibrations lead to the following results.
\begin{thm}\label{mainthm3}
Suppose $\mathfrak M$ is a finite mutually disjoint collection in a {\tt compact} Riemannian manifold $(X, g)$ and every nonempty level $\mathfrak M_k$ satisfy the convex hull condition.
Then for any open neighborhood $U$ of $\underline {\mathfrak M}$,
a new metric $\hat g$ can be constructed by a conformal change of $g$ supported in $U$ such that there exist a family of calibrations
$\{\hat \Phi_{k} \}$ in $(X,\hat g)$ and every nonzero current $T=\sum_{i=1}^{r_k} t_i[[M_i]]$
with $M_i\in\mathfrak M_k$ and 
$t_i\geq0$
is calibrated by $\hat\Phi_k$.
Consequently $T$ is mass-minimizing
in $[T]$
with $\mathrm{\mathbf{M}}(T)=\sum_{i=1}^st_i \mathrm{Vol}_{\hat g}(M_i)$.
\end{thm}
\begin{cor}
%******************
Suppose $\mathfrak M$ is a finite mutually disjoint collection in a {\tt compact} Riemannian manifold $(X, g)$
and
each component represents a nonzero class in the $\mathbb{R}$-homology of $X$.
Then for any open neighborhood $U$ of $\underline {\mathfrak M}$,
 a new metric $\hat g$ can be constructed by a conformal change of $g$ supported in $U$
such that each $\underline {\mathfrak M_k}$ can be tamed in $(X, \hat g)$.
\end{cor}

When $X$ is non-compact, one can obtain better results compared 
with {Theorems \ref{mm}} and {\ref{Xnoncthm}}.
\begin{thm}\label{nc}
Suppose $\mathfrak M$ is a finite mutually disjoint collection in $(X,g)$
and each component represents a nonzero class in the $\mathbb{R}$-homology of $X$. Then
a new metric $\hat g$ can be constructed by a conformal change of $g$
such that every $\underline {\mathfrak M_k}$ can be tamed in $(X,\hat g)$.  
\end{thm}
\begin{thm}\label{nc2}
Suppose $\mathfrak M$ is a neat mutually disjoint collection in $(X,g)$
and each component represents a nonzero class in the $\mathbb{R}$-homology of $X$.
In addition, assume every level of $\mathfrak M$ consist of finite components except the lowest level.
Then a new metric $\hat g$ can be constructed by a conformal change of $g$
such that each $\underline{\mathfrak M_k}$ can be tamed in $(X,\hat g)$.  \end{thm}
\begin{rem}\label{cforccm}
If $(X,g)$ is hermitian with an (almost) complex $J$,
so are
the resulted metrics. 
\end{rem}

As an application, we strengthen Tasaki's ``equivariant" theorem in \cite{Tasaki}.
\begin{thm}[Tasaki] 
Let $K$ be a compact connected Lie transformation group of a manifold $X$ and $M$ be a (connected) compact oriented submanifold in $X$. Assume $M$ is invariant under the action of $K$
and it represents a nonzero $\mathbb R$-homology class of $X$.
Then there exists a $K$-invariant Riemannian metric $g$ on $X$ such that $M$ is mass-minimizing in
homology class with respect to $g$.
\end{thm}

By our method, one can get the following from the proof of {Theorem \ref{eqcase}}.
\begin{thm}\label{eqcase2}
Let $K$ be a compact Lie transformation group of a manifold $X$ and
$M$ be a compact connected oriented submanifold in $X$.
Assume $M$ is invariant under the action of $K$ and the action is orientation preserving.
Then for any $K$-invariant Riemannian metric $g^K$, there exists a $K$-invariant metric $\hat g^K$ conformal to $g^K$ such that $M$ can be calibrated in $(X,\hat g^K)$.
\end{thm}

$K$'s being connected can lead to another generalization.

\begin{thm}\label{eqcase}
Suppose $\frak M$ is a neat mutually disjoint collection
with only the lowest level possibly consisting of infinite components,
and that each component represents a nonzero class in the $\mathbb R$-homology of $X$. 
Let $K$ be a compact connected Lie transformation group of $X$.
Assume $\frak M$ is invariant under the action of $K$.
%(therefore the action is orientation preserving by $K$'s connectedness).
Then for any $K$-invariant Riemannian metric $g^K$, 
there exists a $K$-invariant metric $\hat g^K$ conformal to $g^K$
such that every $\underline {\mathfrak M_k}$ can be tamed in $(X,\hat g^K)$.
\end{thm}

%\begin{rem}
%If $K$ is not connected but it keeps the orientation of $M$ unchanged, then we can still get the same conclusion.
%\end{rem}
 \begin{prooff}
 Without loss of generality, one only needs to consider the case of a single level.
 Since $K$ is compact, there is a {\em Haar}-measure $d\mu$ with $\int_Kd\mu=1$.
 One can use $d\mu$ to average \eqref{abouttoglue} for a $K$-invariant $\Phi$.
 (Note that $\omega^*$ and $\textbf{d}$ are $K$-invariant.)
 Then average the corresponding $\alpha$.
 By \eqref{mgofcc} one can get a $K$-invariant calibration pair 
$(\Phi, \hat g^K)$.
\end{prooff}

\subsection{On Mean Curvature Vector Fields}\label{omcv}
Let us take a short digression about mean curvature vector fields.
By local calibrations, we have the following.

\begin{cor}\label{ctom}
Suppose $M$ is an oriented compact submanifold in $(X, g)$.
Then there exists $\hat g$ conformal to $g$
such that $M$ is minimal in $(X, \hat g)$.
\end{cor}

\begin{rem}
Since either local orientation of a submanifold leads to the same metric by our method
and being minimal is really a local property, the orientability and compactness requirements can be removed.
\end{rem}

What is more, by a direct computation, a concrete relation between mean curvature vector fields
through a conformal change can be given explicitly.

\begin{prop}\label{CCOMC}
Let $M$ be an $m$-dimensional submanifold in $(X,g)$
and $\tilde g=f\cdot g$ where $f$ is a positive function.
Then at a point $p\in M$,
\begin{equation}\label{MC}
f(p)\cdot \tilde H_p=H_p-\frac{m}{2f(p)}\cdot grad_{g,p}^{\perp}(f).
\end{equation}
Here $H$ and $\tilde H$ are mean vector fields of $M$ under $g$ and $\tilde g$ respectively 
and $grad_{g}^{\perp}(\cdot)$ stands for the normal part of $grad_{g}(\cdot)$ along $M$.
\end{prop}

\begin{rem}\label{mcvf}
$M$ can be realized totally geodesic via a
conformal change if and only if M is pointwise totally umbilical.
\end{rem}

\begin{thm}\label{mcvf0}
For a submanifold $M$ (not necessarily oriented or compact) in $(X, g)$ and any (smooth) section $\xi $ of the normal bundle over $M$,
there exists some metric $\tilde g$ conformal to ${g}$ such that %the mean curvature vector field along $M$ with respect to $\tilde {g}$ is $\xi $.
$\tilde H=\xi$.
\end{thm}
\begin{prooff}
Suppose the $\epsilon$-neighborhood $U_\epsilon$ of $M$
for some suitable positive function $\epsilon$ on $M$
can be identified with the normal $\epsilon$-disc bundle
$\mathcal B$ of $M$ via the exponential map restricted to normal directions.
Consider the smooth function $f$ on $\mathcal B$ by
$
f_x(y)=1-\frac{2}{m}<\xi_x-H_g,\ y>_{g_x^\perp}
$
where $x$ is a point of $M$ and $y$ lies in the $\epsilon$-disc fiber through $x$.
Let $F$ be the induced positive (shrink $\epsilon$ if needed) function on $U_\epsilon$.
Take $\check g= F\cdot \hat g$. Since the differential of the identification map
along $M$ is identity, by {\eqref{MC}}
$
%\begin{array}{rcl}
H_{\tilde g}=H_{ g}-\frac{m}{2}\cdot grad^{\perp}_g F
= H_{ g}+\frac{m}{2}\cdot \frac{2}{m}\cdot (\xi_x-H_{g})
= \xi_x.
%\end{array}
$
\end{prooff}
\subsection{Non-Compact $M$}\label{MNC}
Now we consider non-compact submanifolds.
Obviously the mass of a non-compact submanifold usually reaches infinity.
In order to get adapted to this difference, we introduce the following notion. 

\begin{defn}
The (strong) $\bold{Global\ Plateau\ Property}$ of a properly embedded complete non-compact submanifold $M$
means that for any oriented bounded domain $\Omega$ on $M$
with
$d[[\Omega]]\neq 0$
and
$\mathrm{\mathbf{M}}(d[[\Omega]])<\infty$,
$[[\Omega]]$ is a (the unique) area-minimizer among all $m$-dimensional rectifiable currents (see \cite{F} or \cite{FM}) with boundary $d[[\Omega]]$.
\end{defn}

%The following statement asserts a stronger result than being strongly calibrated,
%namely any piece of the submanifold is mass-minimizing among all competitors with the same boundary.

 \begin{thm}\label{GPP}
Suppose $M^m$ is a properly embedded oriented connected non-compact complete
submanifold without boundary in $(X,g)$.
Then for any open neighborhood $U$ of $M$,
a new metric $\hat g$ can be constructed via a conformal change of $g$ supported in $U$
such that $M$ can be strongly calibrated with strong Global Plateau Property in $(X, \hat g)$. 
\end{thm}

\begin{prooff}
Let $U_\epsilon(M)$ be an $\epsilon$-neighborhood of $M$ (for some suitable positive function
$\epsilon$ along $M$) with the bundle structure in \S \ref{BSNS}.
Set $\phi$ $=d(\rho(\textbf{d})\pi_{g}^*(\psi))$
where $d\psi$ is the volume form on $M$,
$\textbf{d}$ and $\rho$ are similar to those in \S\ref{forms}.
Now
$\|\phi\|^*_{g}$ is smooth on $U_{\frac{3}{5}\epsilon}(M)$.
Let $\bar g=(\|\phi\|^*_{g})^2g$.
Then $\|\phi\|^*_{\bar g}=1$ on $U_{\frac{3}{5}\epsilon}(M)$. 
%Therefore we can continue our gluing procedures as before.
Define
$\hat g=\sigma^{\frac{1}{m}}(1+(\textbf{d})^2)\bar g+\alpha(1-\sigma)^{\frac{1}{m}}g$ 
where $\alpha$ is a smooth function with value one on $U_\epsilon^c$
such that $\|\phi\|^*_{\alpha g}\leq 1$ on $X$, and $\sigma$ is the same as in \S\ref{metrics}.
By Lemma \ref{CCGP} and $\mathbf{spt}(\|\phi\|^*_{\hat g}-1)=M$, $\phi$ strongly calibrates $M$ in $(X,\hat g)$. 
Let  $\Omega$ be an oriented bounded domain on $M$
with $\bold{M}(d[[\Omega]])<\infty$.
For any competitor current $K$ (of finite mass)
with $dK=d[[\Omega]]$ we have
$
\mathrm{\mathbf{M}}([[\Omega]])=
[[\Omega]](d(\rho\pi_g^*(\psi)))
=(d [[\Omega]])(\rho\pi_{g}^*(\psi))
=K(\phi)\leq \mathrm{\mathbf{M}}(K).
$
When $K\neq [[\Omega]]$ as functionals over smooth $m$-forms,
the inequality is strict.
\end{prooff}
{\ }

%%%%%%%%%%%%%%%%%%%%%%%%%%%%%%%%
\subsection{Global Plateau Property for Compact $M$}\label{MisC}
How about compact submanifolds? 
By strong Global Plateau Property of a compact submanifold $M$
we mean that for any domain $\Omega$ of $M$
with
$d[[\Omega]]\neq 0$
and
$\mathrm{\mathbf{M}}(d[[\Omega]])<\infty$, $[[\Omega]]$ or (and) $[[\Omega]]-[[M]]$ is (are) the unique (two) area-minimizing rectifiable current(s)
with boundary $d[[\Omega]]$.
      \begin{rem}
      By approximations in $\Omega$ or $\Omega^c$, 
      it can be seen that the requirement $\mathrm{\mathbf{M}}(d[[\Omega]])<\infty$ in the definition of strong Global Plateau Property can be removed.
      \end{rem}
\begin{thm}\label{cptgpp}
Suppose $M$ 
is an oriented compact connected $m$-dimensional submanifold
in $(X,g)$ 
and 
it represents a nonzero class in $H_m(X;\mathbb{R})$. Then
a new metric $\hat g$ can be constructed by a conformal change of $g$
such that  $M$ can be strongly calibrated in $(X, \hat g)$ with strong Global Plateau Property.
\end{thm}
\begin{prooff}
By the compactness of $M$ one can
choose a closed smooth $m$-form $\phi$ 
with $\int_M\phi>0$ and $\|\phi\|^*_g<1$ on $\overline{U_\epsilon}-U_{\frac{4\epsilon}{5}}$.
Suppose one gets $\Phi$ based on $\phi$ as in \S\ref{forms}.
Then there exists $\alpha_\phi$ with value one on $\overline{U_\epsilon}-U_{\frac{4\epsilon}{5}}$
satisfying \eqref{globalcontrolofalpha}.
Denote the resulted metric as in \S\ref{H} by $\hat g$.
Then $\hat g=g$ on $\overline{U_\epsilon}-U_{\frac{4\epsilon}{5}}$.
Assume $\Omega$ is a domain of $M$ with finite $\mathrm{\mathbf{M}}_{\hat g}(d[[\Omega]])$ and
$\bold M_{\hat g}([[\Omega]])< \frac{1}{2}{\text{Vol}}_{\hat g}(M)$. (The proof remains the same for the equality case.)
Suppose $K$ is a rectifiable competitor with
$\mathrm{\mathbf{M}}_{\hat g}(K) < {\text{Vol}}_{\hat g}(M)$
and $dK=d[[\Omega]]$.
Assume also that every connected component of $\mathrm{\bold{spt}}(K)$ touches $\partial \Omega$.

If $\mathrm{\bold{spt}}(K)\subseteq \overline{U_\epsilon}$,
then $[[\Omega]]-K$ belongs to an integer multiple of $[M]$.
Therefore $\mathrm{\mathbf{M}}_{\hat g}([[\Omega]])<\mathrm{\mathbf{M}}_{\hat g}(K)$ by assumption.

If $\mathrm{\bold{spt}}(K)$ is not entirely contained in $\overline{U_\epsilon}$,
a slicing result asserts the existence of very small $0<\mu<\frac{\epsilon}{30}$ such that
the slice $\mathcal S\triangleq <K,\bold{d}, \epsilon-\mu>$ (corresponding to $\bold{d}=\epsilon-\mu$)
is a rectifiable current.
(Note that $\bold d$ is Lipschitz.
This can be achieved by applications of coarea formula. For example, see \cite{FM}.)
A celebrated result in \cite{FF} asserts that $\{\text{integral current } T \text{ with } \mathrm{\bold M} (T)+\mathrm{\bold M} (dT) \leq c \text{ and } \mathrm{\bold{spt}}(T)\subseteq\Gamma\}$
is compact in the weak topology for any positive number $c$ and compact set $\Gamma$.
So there exists an area-minimizing integral current $\tilde K$ with $d \tilde K=[[\Omega]]-\mathcal S$
and $\mathrm{\bold{spt}}(\tilde K)\subseteq \overline{U_{\epsilon-\mu}}$. 
Let $p'$ be a point in $ ({U_{\epsilon-3\mu}}-\overline{U_{\epsilon-4\mu}})\bigcap \mathrm{\bold{spt}}(\tilde K)$.
Set $\overline K_{\Omega,K,\mathcal S,\tilde K}$
to be the restriction $\tilde K\llcorner (U_{\epsilon-\mu}-\overline{U_{\epsilon-6\mu}})$
of $\tilde K$ to ${U_{\epsilon-\mu}-\overline{U_{\epsilon-6\mu}}}$.

\textbf{Claim}: 
There exists some positive $\beta$ independent of  $\Omega$, $K$, $\mathcal S$, $\tilde K$
or $\hat g|_{X-(\overline{U_\epsilon}-U_{\frac{4\epsilon}{5}})}$
such that
$\mathrm{\bold{M}}_{h}(\overline K_{\Omega,K,\mathcal S,\tilde K})>\beta$
for any $\overline K_{\Omega,K,\mathcal S,\tilde K}$ with nonempty support
where $h=g|_{U_{\epsilon-\mu}-\overline{U_{\epsilon-6\mu}}}$.

Let us require in addition at the beginning that
$0<\int_M \phi<\beta$.
Then it follows by the claim that $M$ has strong Global Plateau Property in $(X, \hat g)$.
\end{prooff}

\begin{pf4} The proof follows from Allard's idea.
By Nash's embedding theorem \cite{Nash},
$(\overline{U_\epsilon}-U_{\frac{4\epsilon}{5}},g)$ can be isometrically embedded via a map $f$ into
some Euclidean space $(\mathbb R^s, g_E)$.
%Ignore the subscripts $\Omega,K,\mathcal S,\tilde K$ of $\overline K_{\Omega,K,\mathcal S,\tilde K}$.
$\overline K_{\Omega,K,\mathcal S,\tilde K}$ induces a varifold (introduced by Almgren) $V_{\overline K}$ supported on
${U_{\epsilon-\mu}-\overline{U_{\epsilon-6\mu}}}$.
Since $\overline K$ is a restricted area-minimizing current,
$\delta V_{\overline K}=0$ a.e.
By the compactness of $\overline{U_{\epsilon}}-{U_{\frac{4}{5}\epsilon}}$,
the norm of the corresponding $\delta V_{f_\#\overline K}$ in $\mathbb R^s$ is bounded a.e. by a constant $A$ related to the embedding, namely independent of $\overline K$. 
Note that the density of $V_{f_\#\overline K}$ %(induced by a rectifiable current)
is a.e. at least one on
its support.
Choose a point $a$ on $\mathrm{\bold{spt}}(f_\#\overline K)\bigcap f({U_{\epsilon-3\mu}}-\overline{U_{\epsilon-4\mu}})$ with this density property
in $\mathbb R^s$.
Take $R=\min\{\overline\epsilon,\tau\}$
where $\overline\epsilon$ is some positive number
to ensure
the $\epsilon$-disc bundle structure
over $f(\overline{U_\epsilon}-U_{\frac{4\epsilon}{5}})$ in $\mathbb R^s$ and
$\tau=\mathrm{dist}_{g_E}(f({U_{\epsilon-3\mu}}-\overline{U_{\epsilon-4\mu}}), \partial (f(\overline{U_\epsilon}-U_{\frac{4\epsilon}{5}})))$.
Then our claim follows
by the following monotonicity result in \cite{Allard}.

\begin{thm}[Allard]
Suppose $0\leq A<\infty$, $a\in \text{support of }\|V\|$, $V\in \bold{V}_m(U)$, where $U$ is an open region of $\mathbb R^s$. If $0<R<\text{distance}(a,\mathbb R^s-U)$ and 
$$\|\delta V\|\bold B (a,r)\leq A\|V\|\bold B(a,r)\ \ \ \ \ whenever\ 0<r<R,$$
then $r^{-m}\|V\|\bold B(a,r)\exp Ar$ is nondecreasing in $r$ for $0<r<R$.
\end{thm}
\end{pf4}

\begin{rem}
By our construction, 
if $\mathrm{Vol}_{\hat g}(\Omega)<\frac{1}{2}\mathrm{Vol}_{\hat g}(M)$,
then $[[\Omega]]$ is the unique of least mass among all rectifiable currents with boundary $d[[\Omega]]$.
\end{rem}
% Section 5
{\ }
\\

\section{Several Calibrations}\label{SCal}
%It has been long known since {\em Thom} that, for each $\mathbb Z$-homology class $\alpha$  of an oriented manifold $X$, there exists an integer $N\geq 1$ such that $N\alpha$
%can be represented by some embedded submanifold.
%When $X$ is compact, the sum of betti number is finite.
%Then $N$ can be made universal for elements of some basis of the $\mathbb Z$-homology of every dimension.
In Theorem \ref{mainthm3}, only one calibration
is constructed for each dimension.
Therefore it lacks control on some region in the space of homology classes.
%Note that the sum of two calibrated submanifolds is generally not mass-minimizing.
To conquer this problem
we shall construct
a metric which supports enough calibrations
that we need for each dimension.

When $X^n$ is oriented with betti numbers $b_k<\infty$ for $1\leq k<\frac{1}{2}n$,
by Thom \cite{T} or Corollary II.30 in \cite{Thom}
there exist embedded oriented connected compact  $k$-dimensional submanifolds
$\mathcal L_k\triangleq\{M^k_1,\cdots,M^k_{b_k}\}$
such that $span\{[M^k_i]\}_{i=1}^{b_k}=H_k(X;\mathbb R)$.
By transversality one can assume 
$\bigcup_{1\leq k <\frac{1}{2}n}\mathcal L_k$
is a mutually disjoint collection.
       \begin{thm}\label{1}
       Let $(X^n,g)$ and ${M^k_i}$ be given as above.
       Then there exists $\hat g$ conformal to $g$ such that every nonzero
       $\sum_{i=1}^{b_k} t_i[[M^k_i]]$ 
where $1\leq k<\frac{1}{2}n$, $M^k_i\in\mathcal L_k$ and $t_i\in \mathbb R$ 
%can be strongly calibrated in $(X,\hat g)$.Consequently $\sum_{i=1}^{b_k} t_i[[M^k_i]]$
is the unique mass-minimizing current in
$\sum_{i=1}^{b_k} t_i[M^k_i]$.
\end{thm}
\begin{prooff}
For the sake of simplicity, assume
$\dim H_k(X;\mathbb R)$ $=$ $2$ for some $k< \frac{1}{2}\dim X$
and $\{[M_1],[M_2]\}$ is a basis
where $M_1$ and $M_2$ are disjoint oriented connected compact submanifolds.
Then there exist $k$-forms $\phi_1\text{ and }\phi_2$ on $X$ with
$\int_{M_i}\phi_j=\delta_{i,j}.$
Without loss of generality, assume $\phi_1\equiv 0$ on $U_{\epsilon}(M_2)$ and $\phi_2\equiv 0$ on $U_{\epsilon}(M_1)$.
Note that $\alpha$ can be chosen such that  
$\|\hat \Phi_i\|^*_{\hat g}\leq\frac{1}{2}
%(\text{Vol}_g(M_i))^{\frac{1}{k}}
\text{ on } (U_\epsilon(M_i))^c$
for the constructed metric $\hat g$ and $\hat \Phi_1$ and $\hat \Phi_2$ in \S\ref{ccm}.
%Therefore each $M_i$ is strongly calibrated.
A key observation is that
$
\pm\hat \Phi_1,\ \pm\hat \Phi_2\text{ and }
\pm \hat \Phi_1\pm \hat \Phi_2$
are all calibrations with respect to $\hat g$.

\begin{center}
\includegraphics[scale=0.21]{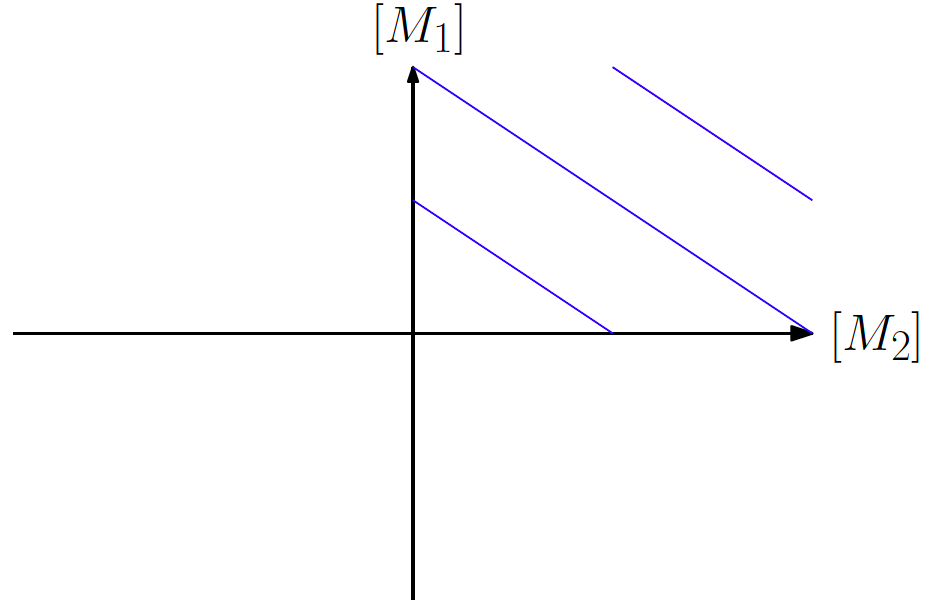}
\end{center}

Obviously any nonzero linear combination of $[[M_1]]$ and $[[M_2]]$
can be calibrated in $(X,\hat g)$.
For examples, those falling into the classes of the closer of the first quadrant can be calibrated by $\hat \Phi_1+\hat\Phi_2$. The uniqueness follows as a result of 
$\mathrm{\bold{spt}}(\|\pm \hat \Phi_i\|^*_{\hat g}-1)= M_i$,
$\mathrm{\bold{spt}}(\|\pm \hat \Phi_1\pm \hat \Phi_2\|^*_{\hat g}-1)=\bigcup M_i$,
the simpleness of $\pm \hat \Phi_i$ along $M_i$
and $\pm \hat \Phi_1\pm \hat \Phi_2$ along $M_1\bigcup M_2$,
and the connectedness of $M_i$.

When $\dim H_k(X;\mathbb R)=s$, $2^s$ such calibrations, each of which has comass norm bounded above by $\frac{1}{s}$ away from some neighborhood of the corresponding submanifold, can be constructed for our purpose.
More generally, for different dimension levels, the above argument combined with the
elimination trick on calibration forms proves the theorem.
\end{prooff}

When the dimension of $X$ is no less than $6$, one can choose $b_k$ smooth $k$-dimensional submanifolds $\{M_i^k\}_{i=1}^{b_k}$ and define
$\mathcal L_k\triangleq \{M^k_{1},\cdots,M^k_{b_k}\}$ for $k=1,2,\cdots, n-3$,
such that
span$\{[M_i^k]\}_{i=1}^{b_k}=H_k(X,\mathbb R)$ and
such that
intersections $\mathcal I$ among $\bigcup^{n-3}_{k=1} \mathcal L_k$ are all transversal.
Note that $\mathcal I$ has a natural stratification structure 
$\cdots\prec \mathcal I_2\prec\mathcal I_1=\bigcup^{n-3}_{k=1} \mathcal L_k$,
where
$\mathcal I_t$ is the set of intersections involving $t$
representatives.

\begin{thm}\label{2}
Let $X^n$ be an oriented manifold
with betti numbers $b_k<\infty$
for $1\leq k\leq n-3$
and $\mathcal L_k$ be given above.
Then there exists a metric $g$
such that 
every nonzero
$\sum_{i=1}^{b_k} t_i[[M^k_i]]$ where $1\leq k\leq n-3$, $M^k_i\in\mathcal L_k$ and $t_i\in \mathbb R$  
is the unique mass-minimizing current in
$\sum_{i=1}^{b_k} t_i[M^k_i]$.
\end{thm}
\begin{prooff}
Choose a metric $g$ on $X$ such that, for any element $S$ of $\mathcal I_t$ ($t\geq 2$),
there exists some ($2\epsilon$-\textbf{cubic}) neighborhood of $S$
with fibers of the bundle structure in \S \ref{BSNS}
split pointwise along $S$ as the Riemannian product of all fibers of $2\epsilon$-neighborhoods of $S$ in
$H_S$ for $H_S \in \mathcal I_{t-1}$ and $S\subseteq H_S$,
and moreover, the horizontal part of $g$ is the pullback of $g|_S$ via the bundle projection.

Let us focus on all the (connected parts of) deepest intersections.
For simplicity,
suppose we have only one connected deepest intersection $\Delta\in\mathcal I_3$, namely
the intersection of three submanifolds $M_1$, $M_2$ and $M_3$.
Assume $2\epsilon$ is universal for $S\in \bigcup_{t\geq 2}\mathcal I_t$ under $g$
in the preceding paragraph.
Denote the volume forms of $M_1$, $M_2$ and $M_3$ by $\omega_1$ $\omega_2$ and $\omega_3$, and
the distance functions to $M_3$ by $\bold{d_3}$. $\bold{d_1}$ and $\bold{d_2}$ are defined similarly.
Since $\omega_i=d\psi_i$ in the $\epsilon$-neighborhood of 
$(M_i\bigcap M_{i+1})\bigcup (M_i\bigcap M_{i+2})$ (subscripts in the sense of mod $3$) in $M_i$,
define $\Psi_i=d(\rho_i\psi_i)$ in the union of $\epsilon$-cubic neighborhoods
of  $M_i\bigcap M_{i+1}$ and $M_i\bigcap M_{i+2}$.
Here we identify the pullback of $\omega_i$ (and $\psi_i$) with themselves, and $\rho_i$ is a smooth increasing function of $\bold{d_i}$ with value zero when $\bold{d_i}\leq \frac{1}{2}\epsilon$ and one when $
\frac{2}{3}\epsilon\leq \bold{d_i}\leq\epsilon$.

\begin{center}
\includegraphics[scale=0.3]{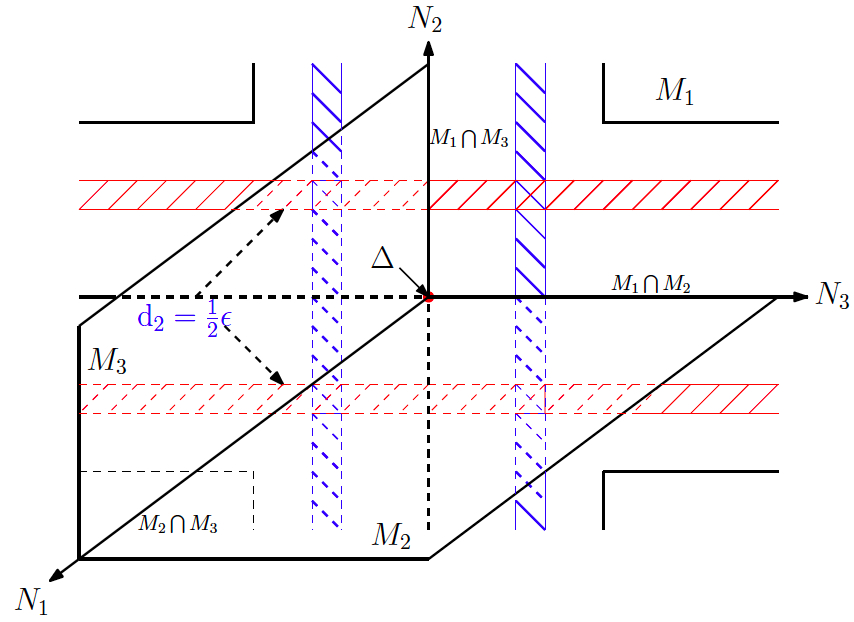}
\end{center}

\textbf{Case One:} Every $M_i$ has the same dimension $k$.
One may want to use the form $\sum\omega_i$.
However it is not well defined on the union $\Xi$ of $\epsilon$-cubic neighborhoods of  $M_1\bigcap M_2$, $M_2\bigcap M_3$ and $M_3\bigcap M_1$.
Let us consider 
\begin{equation*}
\phi_k=\sum \omega_i-\sum \Psi_i=\sum\big[(1-\rho_i)\omega_i-d\rho_i\wedge\psi_i\big] \text{ on } \Xi.
\end{equation*}
Then $\phi_k=\omega_i$ when $\bold{d_{i}}\leq\frac{1}{2}\epsilon$, $\bold{d_{i+1}}\geq\frac{2}{3}\epsilon$
and $\bold{d_{i+2}}\geq\frac{2}{3}\epsilon$.

Since $n-k\geq 3$, $\|\phi_k\|^*=\max\{\|(1-\rho_i)\omega_i-d\rho_i\wedge\psi_i\|^*\}$
pointwise under any metric.
The slash-shadow region and the backslash-shadow region are
the intersections of regions $\frac{1}{2}\epsilon\leq\bold{d_2}\leq\frac{2}{3}\epsilon$
and $\frac{1}{2}\epsilon\leq\bold{d_3}\leq\frac{2}{3}\epsilon$ in $\Xi$ with $M_1$ respectively.
There are three (bunches of) directions $N_i$ on the $\epsilon$-cubic neighborhood $U_\epsilon(\Delta)$ of $\Delta$.
(By the choice of $g$, the existence of $N_i$ cannot be guaranteed when $\bold{d_i}\geq 2\epsilon$.)
Denote the components of $g$ along $N_i$ by $g_i$.

Choose a smooth function $f_i$ of $\bold{d_i}$ on the region
$\bold{d_i}\leq\epsilon$ in $\Xi$ with properties:
(1) $f_i$ equals one when $\bold{d_i}$ $\geq \frac{3}{4}\epsilon$
and is strictly greater than one elsewhere, more precisely
(2) $f_i=1+\bold{d_i}^2$ for $0\leq \bold{d_i}\leq\frac{1}{3}\epsilon$
and 
(3) it has large constant value $C>1$
such that $(\star)$ $C^{-(n-k-2)}\|(1-\rho_i)\omega_i-d\rho_i\wedge\psi_i\|^*_g<1$
on each region
$\Gamma_i:\frac{1}{2}\epsilon\leq\bold{d_i}\leq\frac{2}{3}\epsilon$
in $\Xi$.

For $i=2$, modify $g$ on $U_\epsilon(\Delta)$ in the following way:
$$g_1\rightarrow f_2^{\rho(\bold d_1)}g_1,\quad
g_3\rightarrow f_2^{\rho(\bold d_3)}g_3,\quad
\text{ and } \quad g_2\rightarrow f_2^{-\rho(\bold {d_1})\rho(\bold{d_3})}g_2\quad (\ast)$$
where $\rho$ is a smooth decreasing function on $[0,\epsilon]$ with value one on $[0,\frac{4}{5}\epsilon]$
and zero on $[\frac{5}{6}\epsilon,\epsilon]$.
This modification extends trivially to
the union $\Upsilon$ of $\frac{4}{5}\epsilon$-cubic neighborhoods of  $M_1\bigcap M_2$, $M_2\bigcap M_3$ and $M_3\bigcap M_1$,
and moreover, it preserves the volume form of $M_i\bigcap\Upsilon$.

Multiply each of $g_1$, $g_2$ and $g_3$ by the product of corresponding conformal factors in $(\ast)$ for $g$ and $i=1,2,3$.
Denote the new metric on $\Upsilon$ by $\check g$.
The last property of $f_i$ and
$\|\omega_i-\Psi_i\|^*_{\check g}\leq
C^{-(n-k-2)}\|\omega_i-\Psi_i\|^*_{ g}$
on $\Gamma_i$
in fact shows that nonzero
$\sum_{i, n_i=\pm 1,0}n_i(\omega_i-\Psi_i)$
become calibrations when restricted to $(\Upsilon, \check g)$.
Furthermore they naturally extend to 
calibration pairs on
the union $\Theta$ of $\Upsilon$ and the $\frac{1}{2}\epsilon$-neighborhood of $\bigcup M_i$
with respect to the trivial extension $\overline g$ of $\check g$.
\\{\ }
 
\textbf{Case Two:}  $M_2$ and $M_3$ are of dimension $k$, but $M_1$ has a different dimension $m$.
(Similar for the case with mutually different dimensions.)
Consider potential calibrations
$ \pm(\omega_2-\Psi_2), \pm(\omega_3-\Psi_3),\pm(\omega_2-\Psi_2)\pm(\omega_3-\Psi_3),
\text{ and }\pm(\omega_1-\Psi_1)$
on $\Xi$.
By the same procedure (with different weights in $(\star)$ and $(\ast)$),
one can get similar calibration paris
in some neighborhood of $\bigcup M_i$.
\\{\ }

Clearly the idea works for general case (with modifications in $(\star)$ and $(\ast)$).
%Now extend local calibrations to global ones.
Following the above procedure around all connected parts of deepest intersection, one can
extend these finite preferred local calibration pairs to global ones sharing a common metric.
Multiply the metric by a smooth function which is one along $\bigcup_{k=1}^{n-3}\bigcup_{M\in\mathcal L_k}M$
and strictly greater than one elsewhere.
Name the final metric
$\hat g$.
By our construction, every nonzero $\sum_{t_i\in\mathbb R,M_i\in\mathcal L_k} t_i[[M_i]]$
with $1\leq k\leq n-3$
can be calibrated in $(X,\hat g)$. The uniqueness of such mass-minimizing current in its current homology class follows
similarly as in the proof of Theorem \ref{2}.
Note that
for any point $p\in M^k_i-\mathcal I_2$
the oriented unit (volume) $k$-vector of $\wedge^k{T_pM_i^k}$
is the unique among all unit $k$-vectors in $\wedge^kT_pX$
such that
its pairing with the corresponding calibration of $M_i^k$ produces value one.
\end{prooff}
  
\begin{rem}
For arbitrary metric $h$,
one can still use the (local) bundle structures induce by $g$.
For example, consider Case One above.
Modify $h$ first as follows.
Let $a_i$ denote $\|\omega_i\|_{h}^*$ in $U_{2\epsilon}(\Delta)$.
Then multiply $h$ along $N_i$ by $(\frac{a_{i+1}a_{i+2}}{a_i})^\frac{1}{n-k}$ for $i=1,2,3$.
Now $\|\omega_i\|_{\bar h}^*=1$ under the new metric $\bar h$ .
Based on this, one can get a metric $\tilde h$ 
on the union of $\epsilon$-neighborhood $\Sigma_i$ of $M_i$
(obtained by conformal changes of $h$
along ``split horizontal directions"
on $\Sigma_i$)
such that $\|\omega_i\|_{\tilde h}^*\equiv 1$
on $\Sigma_i$.
(Extension can be made through $\{\epsilon<\bold{d_i}<2\epsilon\}\bigcap U_{2\epsilon}(\Delta)$.)
Then one can apply the same construction based on $\tilde h$.
\end{rem}
{\ }
   \begin{rem}
     The requirement of codimension no less than $3$ is crucial in our proof. 
     When $n=4$ or $5$, Theorem \ref{2} can be improved to include the level of codimension $2$ by \cite{Z4}.
     \end{rem}
{\ }
\\
\textbf{Open Question:}
Usually we cannot have such existence result when $k$ can be $n-1$.
Therefore it may be interesting to study whether the same result holds for $1\leq k\leq n-2$.  
{\ }

%%%%%%%%%%%%%%%%%%%%%%%%%%%%%%%
%\clearpage
\appendix
{\ }

\section{\sc A Generalized Gauss Lemma}\label{C}

Since any one-dimensional foliation is always locally orientable, one can define a local \textbf{length flow} with respect to a metric according to a choice of orientation.
An observation of Sullivan \cite{S2} and a special case in Harvey and Lawson \cite{HL1} are the following. 

\begin{thm}[Sullivan]
A one-dimensional flow is geodesible if and only if there is a transverse field of codimension one
planes invariant under the length flow.
\end{thm}
\begin{thm}[Harvey and Lawson]
A one-dimensional foliation $\Gamma$ is geodesic if and only if its perpendicular plane field $\mathscr P$ is invariant under the length flow.
\end{thm}
%\begin{center}
%\includegraphics[scale=0.45]{f12.jpg}
%\end{center}
\begin{prooff}
Locally, denote the oriented unit tangent vector field by $V$.
For any local (nowhere zero) smooth section $N$ of $\mathscr P \cong TX/T\Gamma$, we have
$0=N<V,V>=2<\nabla_NV,V>.$
Furthermore,
\begin{equation}\label{apc}
\begin{split}
0&=V<V,N>\
=\ <\nabla_VV,N>+<V,\nabla_VN>\\
&=<\nabla_VV,N>+<V,\nabla_NV+\mathfrak{L}_VN>\\
&=<\nabla_VV,N>+<V,\mathfrak{L}_VN>.
\end{split}
\end{equation}
Hence
$\Gamma$ is geodesic if and only if $\mathscr P$ is preserved by the local length flow.
\end{prooff}

Since foliations that we encounter, e.g. in the proof of {Propersition \ref{wps}}, are all locally integrable,
we focus on this special situation.

\begin{cor}[Generalized Gauss Lemma]\label{ggl}
Suppose $\Gamma$ is a one-dimensional foliation % of an open ball $B^n$ endowed with a metric $g$.
and its perpendicular plane field $\mathscr P$ is locally integrable.
Then $\Gamma$ is geodesic if and only if 
local integral pieces of $\mathscr P$
are preserved by the local length flow along $\Gamma$.
\end{cor}

%\begin{rem}\label{rkggl}
%Suppose that the local perpendicular plane field $\mathscr P$of a geodesic foliation $\Gamma$ through $p$ is integrable.
%Then $\mathscr P$ is integrable around $p$.  
%\end{rem}
\begin{rem}
The appellation ``Generalized Gauss Lemma" is due to the fact that
one can derive Gauss Lemma from the corollary and \eqref{apc}.
\end{rem}

\begin{cor}\label{corggl}
Suppose $\Gamma$ is a one-dimensional geodesic foliation %of an open ball $(B^n,g)$ and 
and its perpendicular plane field $\mathscr P$ is locally integrable.
Let $g=g^\Gamma\oplus g^\perp$ be the metric decomposition of $g$ along $\Gamma$ and $\mathscr P$.
Assume $\hat g=g^\Gamma\oplus {\hat g}^\perp$ is a smooth metric by replacing $g^\perp$ by $\hat g ^\perp$.
Then $\Gamma$ is geodesic as well with respect to $\hat g$.
\end{cor}
\begin{prooff}
Since $g^\Gamma$ and $\Gamma$'s perpendicular plane field are unchanged,
$\Gamma$ and $\hat g$ satisfy the conditions in Corollary\ \ref{ggl}.
So the conclusion follows.
\end{prooff}

{\ }
\\
{\  }


\begin{bibdiv}
\begin{biblist}




\bib{Allard}{article}{
    author={Allard, William K.}
    title={On the First Variation of a Varifold},
    journal={Ann. Math. },
    volume={95},
    date={1972},
    pages={417--491},
}

\bib{F}{book}{
    author={Federer, Herbert},
    title={Geometric Measure Theory},
    place={Springer-Verlag, New York},
    date={1969},
}

\bib{FF}{article}{
    author={Federer, Herbert},
    author={Fleming, Wendell H.},
    title={Normal and integral currents},
    journal={Ann. Math. },
    volume={72},
    date={1960},
    pages={458--520},
}

%\bibitem{Gluck} H. Gluck, {\em Open letter on geodesible flows.}


\bib{Gluck}{article}{
    author={Gluck, Herman},
    title={Open letter on geodesible flows.}
}

%\bibitem{HL1} R. Harvey and H. B. Lawson, Jr.,\ \ {\em Calibrated Foliations},\  Amer. J. Math. (1982) Vol. 104:\ 607-633.

\bib{HL1}{article}{
    author={Harvey, Reese},
    author={Lawson, Blaine H.},
    title={Calibrated Foliations},
    journal={Amer. J. Math.},
    volume={104},
    date={1982},
    pages={607--633},
}
%\bibitem{HL2} R. Harvey and H. B. Lawson, Jr.,\ \ {\em Calibrated geometries},\ Acta Math. (1982) Vol. 148:\ 47-157.
\bib{HL2}{article}{
    author={Harvey, Reese},
    author={Lawson, Blaine H.},
    title={Calibrated geometries},
    journal={Acta Math.},
    volume={148},
    date={1982},
    pages={47--157},
}

%\bibitem{FM} F. Morgan,\ \ {Geometric Measure Theory: A Beginner's Guide}, {Academic Press}; 4th edition, 2008.

\bib{FM}{book}{
    author={Morgan, Frank},
    title={Geometric Measure Theory: A Beginner's Guide},
    place={Academic Press, 4th edition},
    date={2008},
}

%\bibitem{Nash} J. Nash,\ \ {\em The imbedding problem for Riemannian manifolds},\ Ann. of Math. (1956) Vol. 63:\ 20-63.

\bib{Nash}{article}{
    author={Nash, John},
    title={The imbedding problem for Riemannian manifolds},
    journal={Ann. of Math.},
    volume={63},
    date={1956},
    pages={20--63},
}


%\bibitem{S2} D. Sullivan,\ \ {\em A foliation by geodesics is characterized by having no ``tangent homologies"}, J. Pure and Applied Algebra (1978) Vol. 13:\ 101-104.

\bib{S2}{article}{
    author={Sullivan, Dennis},
    title={A foliation by geodesics is characterized by having no ``tangent homologies"},
    journal={J. Pure and Applied Algebra},
    volume={13},
    date={1978},
    pages={101--104},
}


%\bibitem{S3} D. Sullivan,\ \ {\em A homological characterization of foliations consisting of minimal surfaces}, Comm. Math. Helv. (1979) Vol. 54:\ 218-223.
\bib{S3}{article}{
    author={Sullivan, Dennis},
    title={A homological characterization of foliations consisting of minimal surfaces},
    journal={Comm. Math. Helv.},
    volume={54},
    date={1979},
    pages={218--223},
}


%\bibitem{Tasaki} H. Tasaki,\ \ {\em mass minimizing submanifolds with respect to some Riemannian metrics},\ J. Math. Soc. Japan (1993) Vol. 45:\ 77-87.
\bib{Tasaki}{article}{
    author={Tasaki, Hiroyuki},
    title={mass minimizing submanifolds with respect to some Riemannian metrics},
    journal={J. Math. Soc. Japan},
    volume={45},
    date={1993},
    pages={77--87},
}

%\bibitem{T} R. Thom,\ \ {\em Quelques propri$\acute{e}$t$\acute{e}$s globales des variŽtŽs diff$\acute{e}$rentiables,} Comment. Math. Helv. (1954) Vol. 28:\ 17-86.
\bib{T}{article}{
    author={Thom, Ren\'e},
    title={Quelques propri$\acute{e}$t$\acute{e}$s globales des variŽtŽs diff$\acute{e}$rentiables},
    journal={Comm. Math. Helv.},
    volume={28},
    date={1954},
    pages={17--86},
}

%\bibitem{Thom} R. Thom,\ \ {Some ``global" properties of differentiable manifolds} (Translated by V. O. Manturov with M. M. PostnikovÕs comments (1958)), pp. 131-209, in
%{Topoligical Library I: Cobordisms and Their Applications}, Editors(s): Novikov,
%World Scientific, Singapore, 2007.

\bib{Thom}{book}{
   author={Thom, Ren\'e},
    title={Some ``global" properties of differentiable manifolds (Translated by V. O. Manturov with M. M. PostnikovÕs comments (1958))},
   place={pp. 131-209, in {Topoligical Library I: Cobordisms and Their Applications}, Editors: Novikov,
World Scientific, Singapore},
   date={2007},
}

%\bibitem{Z4} Y. Zhang,\ \ {\em On extending calibrations,} (part of ``On Calibrated Geometry I: Gluing Techniques").

\bib{Z4}{article}{
   author={Zhang, Yongsheng},
   title={On extending calibrations (part of ``On Calibrated Geometry I: Gluing Techniques")}
   }

\end{biblist}
\end{bibdiv}
\end{document}